\documentclass[12pt]{article}

\usepackage{amssymb}
\usepackage{amsfonts}
\usepackage{amsmath}
\usepackage{theorem}
\usepackage{amscd}
\input xy
\xyoption{all}

\theoremstyle{plain}
\newtheorem{thm}{Theorem}[section]

\newtheorem{lem}[thm]{Lemma}
\newtheorem{coro}[thm]{Corollary}

\newtheorem{prop}[thm]{Proposition}

\theoremstyle{definition} \theoremstyle{definition}
\newtheorem{defn}[thm]{Definition}

\newtheorem{rem}[thm]{Remark}           

\theoremstyle{remark}

\newenvironment{proof}{\noindent {\em Proof.}}{$\Box$ \newline}

\newtheorem{step}{Step}

\def\vertdot{|\, .\,|}

\newcommand{\C}{\mathbb C}
\newcommand{\Q}{\mathbb Q}
\newcommand{\Z}{\mathbb Z}
\newcommand{\R}{\mathbb R}
\newcommand{\Ind}{{\rm Ind}\,}

\def\M{{\mathsf M}}

\newcommand{\Hom}{{\rm Hom}}

\newcommand{\PGL}{{\rm PGL}}
\newcommand{\Sp}{{\rm Sp}}
\newcommand{\GSO}{{\rm GSO}}

\def\A{{\mathbb A}}
\def\C{{\mathbb C}}

\def\Q{{\mathbb Q}}

\def\Z{{\mathbb Z}}
\def\C{{\mathbb C}}
\def\N{{\mathbb N}}

\def\e{{\mathbf e}}
\def\O{{\mathcal O}}

\def\GO{{\rm GO}}
\def\SO{{\rm SO}}
\def\GL{{\rm GL}}
\def\PGL{{\rm PGL}}

\def\Ind{{\rm Ind}}

\def\W{{\mathcal W}}

\def\zZ{{\mathcal Z}}

\def\tr{{\rm tr}}

\def\GSp{{\rm GSp}}

\makeatother \makeatletter
\begin{document}
\title{Spinor $L$-Functions, Theta Correspondence, and Bessel Coefficients}
\author{Ramin Takloo-Bighash}
\maketitle

The purpose of this note is to establish the entireness of spinor
$L$-function of certain automorphic cuspidal representations of
the group similitude symplectic group of order four over the
rational numbers. Our study of spinor $L$-function is based on an
integral representation which works only for generic
representations. For this reason, while methods of this papers do
not directly apply to the most interesting case of interest, ie.
Siegel modular forms of genus two, they show what ``should" be
true for holomorphic forms; after all, generic forms are expected
to be in a certain sense typical. These integrals which were
introduced by M. Novodvorsky in the Corvallis conference \cite
{Novodvorsky} serve as one of the few available integral
representations for the Spinor $L$-function of $\GSp(4)$. Some of
the details missing in Novodvorsky's original paper have been
reproduced in Daniel Bump's survey article \cite {Bump}. Further
details have been supplied by \cite{Takloo-Bighash}. David Soudry
has generalized the integrals considered in Novodvorsky's paper to
orthogonal groups of arbitrary odd degree.

In light of the results of \cite{Takloo-Bighash}, it is sufficient
to study the integral of Novodvorsky at the archimedean place.
Archimedean computations are often forbidding, and unless one
expects major simplifications due to the nature of the parameters,
the resulting integrals are often quite hard to manage. In our
case of interest, the work of Moriyama \cite{Moriyama} benefits
from exactly such simplifications when he treats the case of
cuspidal representations with archimedean components in the
generic (limit of) discrete series. In this work, we concentrate
on those archimedean representations for which direct computations
have yielded very little. For this reason, our methods are a bit
indirect, in fact somewhat more indirect that what at first seems
necessary. Our method is based on the theta correspondence. First
we observe that Novodvorsky's integral is in fact a split Bessel
functional. Then we pull back the Bessel functional via the theta
correspondence for the dual reductive pair $(\GO(2, 2), \GSp(4))$,
and prove that the resulting functional on $\GO(2,2)$ is Eulerian.
On the other hand, one can prove that the integral of Novodvorsky
itself is Eulerian, with an Euler product involving the Whittaker
functions. Next obvious step is to pull back the Whittaker
function via the theta correspondence. Now we have obtained two
different Euler product expansions which represent the same
object, but do not look the same. Then one uses the standard
technique of twisting with highly ramified characters to isolate
the archimedean place to obtain an identity expressing the local
Novodvorsky integral at the archimedean place in terms of an
expression which does not go through the local Whittaker
functions. The advantage of using this expression is that, first
it avoids Whittaker function, so it is effectively more
elementary, and second one can devise a two complex variable zeta
function to study its analytic properties. This identity, at
first, is established only for those representations which appear
as archimedean components of global theta lifts from $\GO(2,2)$.
Then one uses various density arguments to extend the identity to
other cases. The next natural step is to examine the identity to
see what representations of the archimedean group have been
covered. I suspect that at least all unramified tempered
representations are included.

As mentioned above, the main contribution of this work, if any, is
the archimedean analysis. Some of the results of this paper,
especially in the case of discrete series representations, were
announced in \cite{Takloo-Bighash3}. As stated above, the
appearance of \cite{Moriyama} has made our results for discrete
series representations obsolete; Moriyama has obtained better and
more explicit results for generic (limits of) discrete series, and
some other representations, using more direct methods. Also we
have recently learned that Asgari and Shahidi are preparing a
manuscript which contains, among other things, the lifting of
generic automorphic forms from spinor groups to general linear
groups; if established, our results would be trivialized, as
$\GSp(4)$ is nothing but ${\rm GSin}_5$. It appears that Brooks
Roberts has used methods very similar to ours in \cite{Roberts2}
to study various non-archimedean questions. It turns out that both
of us were influenced by Masaaki Furusawa, and communication with
him and Shalika was our common source of inspiration. I learned
about Bessel functionals and theta correspondence from J. A.
Shalika while a graduate student at Johns Hopkins. Here we thank
Shalika for continued support and encouragement over the past few
years. Most of preliminary computations that led to the writing of
this paper were also performed at Hopkins under his supervision.
The author has benefited from conversations with Freydoon Shahidi,
Mahdi Asgari, Jeffrey Adams, Akshay Venkatesh, Peter Sarnak, and
Brooks Roberts for answering many questions regarding his papers.
Also the author wishes to thank the Clay Mathematical Institute
for partial support of the project.
\tableofcontents
\subsection*{Notation} In this paper, the group $\GSp(4)$ over an
arbitrary field $K$ is the group of all matrices $g \in \GL_4(K)$
that satisfy the following equation for some scalar $\nu (g) \in
K$:
\begin{equation*}
  ^t g J g
= \nu (g) J,
\end{equation*}
where $J=\begin{pmatrix} &&1 \\ &&&1 \\ -1 \\ & -1 \end{pmatrix}$.
It is a standard fact that $G=\GSp(4)$ is a reductive group. The
map $(F^{\times})^3 \longrightarrow  G$, given by
\begin{equation*}
(a, b, \lambda) \mapsto \text{diag} (a, b, \lambda a^{-1}, \lambda
b^{-1})
\end{equation*}
gives a parameterization of a maximal torus $T$ in $G$. Let
$\chi_1$, $\chi_2$ and $\chi_3$ be quasi-characters of
$F^{\times}$. We define the character $\chi_1 \otimes \chi_2
\otimes \chi_3$ of $T$ by
\begin{equation*}
(\chi_1 \otimes \chi_2 \otimes \chi_3) ( \text{diag} (a, b,
\lambda a^{-1}, \lambda b^{-1})) = \chi_1(a) \chi_2(b)
\chi_3(\lambda).
\end{equation*}
The Weyl group is a dihedral group of order eight. We have three
standard parabolic subgroups: The Borel subgroup $B$, The Siegel
subgroup $P$, and the Klingen subgroup $Q$ with the following Levi
decompositions:
\begin{equation*}
B= \left\{
\begin{pmatrix} a \\ & b \\ && a^{-1}\lambda \\ &&&b^{-1} \lambda
\end{pmatrix}
\begin{pmatrix} 1 & x \\ & 1 \\ &&1& \\ &&-x&1 \end{pmatrix}
\begin{pmatrix} 1&&s&r \\ &1&r&t \\ &&1 \\ &&&1 \end{pmatrix}
\right\},
\end{equation*}
\begin{equation*}
P= \left\{
\begin{pmatrix} g \\ & \alpha \,^T g ^{-1} \end{pmatrix}
\begin{pmatrix} 1&&s&r \\ &1&r&t \\ &&1 \\ &&&1 \end{pmatrix}; g \in \GL(2) \right\},
\end{equation*}
and finally $Q$ is the maximal parabolic subgroup with non-abelian
unipotent radical associated to the long simple root. Over a local
field, we will use the notation $\chi_1 \times \chi_2 \rtimes
\chi_3$ for the parabolically induced representation from the
minimal parabolic subgroup, by the character  $\chi_1 \otimes
\chi_2 \otimes \chi_3$. If $\pi$ is a smooth representation of
$\GL(2)$, and $\chi$ a quasi-character of $F^{\times}$, then $\pi
\rtimes \chi$ (respectively $\chi \rtimes \pi$) is the
parabolically induced representation from the Levi subgroup of the
Siegel (resp. Klingen) parabolic subgroup. We define a character
of the unipotent radical $N(B)$ of the Borel subgroup by the
following:
\begin{equation*}
\theta (\begin{pmatrix} 1 & x \\ & 1 \\ &&1& \\ &&-x&1
\end{pmatrix}
\begin{pmatrix} 1&&s&r \\ &1&r&t \\ &&1 \\ &&&1 \end{pmatrix}) = \psi(x+t).
\end{equation*}
We call an irreducible representation $(\Pi, V_\Pi)$ of $\GSp(4)$
over a local field generic, if there is a functional $\lambda_\Pi$
on $V_\Pi$ such that
\begin{equation*}
\lambda_\Pi(\Pi(n)v) = \theta(n) v,
\end{equation*}
for all $v \in V_\Pi$ and $n \in N(B)$. If such a functional
exists, it is unique up to a constant \cite {Shalika}. Freydoon
Shahidi has given canonical constructions of these functionals in
\cite {Shahidi} for representations induced from generic
representations. We define Whittaker functions on $G \times V_\Pi$
by
\begin{equation*}
W(\Pi, v, g) = \lambda_\Pi(\Pi(g)v).
\end{equation*}
When there is no danger of confusion, after fixing $v$ and
suppressing $\Pi$, we write $W(g)$ instead of $W(\Pi, v, g)$. For
any representation $\pi$, we will denote by $\omega_\pi$ the
central character of $\pi$.
\section{The work of Jacquet and Langlands}
In this section, we concentrate on the simpler group $\GL(2)$.
This section serves as motivation for later sections which contain
the main results of the paper. Our exposition is heavily based on
\cite{Gelbart-Shahidi}, to the point of copying, especially pages
5-19. For the sake of familiarity and simplicity, we work over
$\Q$.

Let $G= \GL(2)$. Suppose $\chi$ is a unitary character of
$\A^\times$. By a $\chi$-cusp form $\varphi$ on $GL(2)$, we mean
an $L^2(Z_\A G(\Q)\backslash G(\A))$ function satisfying
\begin{equation*}
\varphi( \begin{pmatrix} a \\ & a \end{pmatrix} g ) = \chi(a)
\varphi(g),
\end{equation*}
and
\begin{equation*}
\int_{\Q \backslash \A} \varphi ( \begin{pmatrix} 1 & x \\ & 1
\end{pmatrix} g ) \, dx = 0,
\end{equation*}
for almost all $g \in G(\A)$. It is clear that if $\varphi$ is a
$\chi$-cusp form, and $g \in G(\A)$, then the function $g.
\varphi$ on $G(\A)$ defined by
\begin{equation*}
g.\varphi(h) = \varphi(h g),
\end{equation*}
is again a $\chi$-cusp form. This defines a representation of
$G(\A)$ on the vector space of $\chi$-cusp forms $L_0^2(\chi)$. It
is a fundamental fact that $L_0^2(\chi)$ is a discrete direct sum
of irreducible subspaces each of which appears with multiplicity
one. An irreducible representation $\pi$ of $\GL_2(\A)$ which is
realized as an irreducible subspace $H_\pi$ of $L_0^2(\chi)$ is
called a {\em cuspidal automorphic representation}.

Suppose $\pi$ is an irreducible cuspidal automorphic
representation of $\GL_2(\A)$, and $\varphi \in H_\pi$. We
introduce a {\em global zeta integral}
\begin{equation}\label{eq:Zeta}
\zZ(\varphi, s) = \int_{\Q^\times \backslash \A^\times} \varphi
(\begin{pmatrix} a \\ & 1 \end{pmatrix} ) \vert a \vert_\A^{s - {1
\over 2}}\, d^\times a.
\end{equation}
If $\varphi$ is ``nice enough", $\zZ(\varphi, s)$ defines an
entire function in $\C$. Also, the zeta function $\zZ(\varphi, s)$
satisfies a functional equation:
\begin{equation*}
\zZ(\varphi, s) = \tilde{\zZ}(\varphi^w, 1-s),
\end{equation*}
 where $w = \begin{pmatrix} & 1 \\
-1 \end{pmatrix}$, and $\varphi^w(g) = \varphi(gw)$. Also
\begin{equation}
\tilde{\zZ}(\varphi, s) = \int_{\Q^\times \backslash \A^\times}
\varphi( \begin{pmatrix} a \\ & 1 \end{pmatrix}) \vert a \vert^{s
- {1 \over 2}}\chi^{-1}(a)\, d^\times a,
\end{equation}
with $\chi$ the central character of $\pi$.

The problem is to relate the function $\zZ(\varphi, s)$ to an
automorphic L-function $L(s, \pi, r)$ for some representation $r$
of $^L G = \GL_2(\C)$. For this purpose, we start by the Fourier
expansion of $\varphi$
\begin{equation}\label{eq:Fourier}
\varphi(g) = \sum_{\xi \in \Q^\times} W_\varphi^\psi(
\begin{pmatrix} \xi \\ & 1 \end{pmatrix} g).
\end{equation}
Here
\begin{equation}
W_\varphi^\psi(g) = \int_{\Q \backslash \A} \varphi
(\begin{pmatrix} 1 & x \\ & 1 \end{pmatrix} g) \overline{\psi(x)}
\, dx,
\end{equation}
where $\psi$ is a non-trivial character of $\Q \backslash \A$. It
follows from (\ref{eq:Zeta}) and (\ref{eq:Fourier}) that
\begin{equation*}
\zZ(\varphi, s) = \int_{\A^\times} W_\varphi^\psi( \begin{pmatrix}
a \\ & 1 \end{pmatrix}) \vert a \vert^{s - {1 \over 2}}\, d^\times
a,
\end{equation*}
for $\Re s$ large enough.

Now we recall some of the properties of the {\em Whittaker
functions} $W_\varphi^\psi$. From now on we suppress $\psi$. We
assume that $\varphi$ is right K-finite. In this situation,
$W_\varphi$ is rapidly decreasing at infinity and satisfies
\begin{equation}
W_\varphi( \begin{pmatrix} 1 & x \\ &1 \end{pmatrix} g ) = \psi(x)
W_\varphi(g),
\end{equation}
for all $x \in \A$. The space of all such $W_\varphi$ provides the
$\psi$-Whittaker model of $\pi$. It is known that such a model is
unique (\cite{Jacquet-Langlands}, or \cite{Shalika}), and it is
equal to the restricted tensor product of local Whittaker models
$\mathcal{W}(\pi_p, \psi_p)$, where $\pi = \otimes_p \pi_p$ and
$\psi = \prod_p \psi_p$. In particular, we can assume that
\begin{equation}\label{eq:Whittaker}
W_\varphi(g) = \prod_p W_p(g_p),
\end{equation}
where each $W_p \in \mathcal{W}(\pi_p, \psi_p)$ and for almost all
finite $p$, $W_p$ is {\em unramified}, i.e. $W_p(k) = 1$ for $k
\in K_p = \GL_2(\Z_p)$.

Finally we obtain for $\Re s$ large
\begin{equation}
\zZ(\varphi, s) = \prod_p \zZ(W_p, s),
\end{equation}
where
\begin{equation}
\zZ(W_p, s) = \int_{\Q_p^\times} W_p( \begin{pmatrix} a \\ & 1
\end{pmatrix})\vert a \vert^{s - {1 \over 2}}\, d^\times a.
\end{equation}
First, we collect some of the properties of the local zeta
functions $\zZ(W, s)$. The fundamental fact for $p < \infty$ is
the following: There are a finite number of finite functions $c_1,
\dots, c_N$ on $\Q_p^\times$, depending only on $\pi_p$, such that
for every $W \in \mathcal{W}(\pi_p, \psi_p)$, there are
Schwartz-Bruhat functions $\Phi_1, \dots, \Phi_N$ on $\Q_p$
satisfying
\begin{equation}\label{eq:asymptotic}
W(\begin{pmatrix} a \\ & 1 \end{pmatrix} ) = \sum_{i=1}^N c_i(a)
\Phi_i(a).
\end{equation}
Here, a finite function is a function whose space of right
translates by $\Q_p^\times$ is finite dimensional; finite
functions on $\Q_p^\times$ are thus characters, integer powers of
the valuation function, or products and linear combinations
thereof. Taking the asymptotic expansion just mentioned for
granted, we obtain from Tate's thesis that the integral defining
$\zZ(W, s)$ converges for $\Re s$ large (independent of $s$), and
in the domain of convergence is equal to a rational function of
$p^{-s}$. In particular, the integral has a meromorphic
continuation to all of $\C$. Furthermore, the family of rational
functions $\{ \zZ(W, s) \, \vert \, W \in \mathcal{W}(\pi_p,
\psi_p)\}$ admits a common denominator, i.e. a polynomial $P$ such
that $P(p^{-s})\zZ(W, s) \in \C [p^{-s}, p^s ]$, for all $W$.
Also, there exists a $W^*$ in $\mathcal{W}(\pi_p, \psi_p)$ with
the property that $\zZ(W^*, s) = 1$. The analogous result in the
archimedean situation is that there is a $W^*$ such that $\zZ(W^*,
s)$ has no poles or zeroes.

As in Tate's thesis, we also need to establish the functional
equation and perform the local unramified computations. The
functional equation asserts that there exists a meromorphic
function $\gamma(\pi_p, \psi_p, s)$ (rational function in $p^{-s}$
when $p < \infty$) such that
\begin{equation}
\tilde{\zZ}(W^w, 1-s) = \gamma(\pi_p, \psi_p, s) \zZ(W, s).
\end{equation}
Here $W^w(g) = W(gw)$ and
\begin{equation*}
\tilde{\zZ}(W, s) = \int_{\Q_p^\times} W( \begin{pmatrix} a \\ & 1
\end{pmatrix}) \vert a \vert^{s - {1 \over 2}} \chi_p^{-1}(a)\,
d^\times a,
\end{equation*}
with $\chi_p$ the central character of $\pi_p$. The (non-trivial)
proof of this equation follows from the fact that the integrals
$\zZ$ and $\tilde{\zZ}$ define functionals (depending on $s$)
satisfying a certain invariance property. Next one proves that the
space of such functionals is at most one dimensional, implying
that $\zZ$ and $\tilde{\zZ}$ must be proportional. The factor
$\gamma$ is simply the factor of proportionality.

We recall $\pi_p$ is called unramified when
\begin{equation*}
\pi_p = {\rm Ind}(\mu_1 \otimes \mu_2),
\end{equation*}
with $\mu_1$ and $\mu_2$ unramified characters of $\Q_p^\times$.
We also suppose that $W^0$ is the unique $K_p$-invariant function
in $\mathcal{W}(\pi_p, \psi_p)$. We also suppose that $\psi_p$ is
unramified. Then a direct calculation, using the results of
\cite{Casselman-Shalika} for example, shows that
\begin{equation*}
\zZ(W^0, s) = \frac{1}{(1 - \mu_1(\varpi_p) p^{-s})(1-
\mu_2(\varpi_p)p^{-s})},
\end{equation*}
where $\varpi_p$ is the local uniformizer at $p$. Next the
conjugacy class in $^L G = \GL_2(\C)$ canonically associated with
$\pi_p$ is $t_p = \begin{pmatrix} \mu_1(\varpi_p) \\ &
\mu_2(\varpi_p) \end{pmatrix}$. In particular, we have
\begin{equation}\label{eq:gcd}
\zZ(W^0, s) = L_p(s, \pi, r),
\end{equation}
with $r$ the standard two dimensional representation of
$\GL_2(\C)$.

After this preparation, we can prove the conjecture of Langlands
for $L(s, \pi, r)$ with $r$ as above. For simplicity, we write
$L(s, \pi)$ instead of $L(s, \pi, r)$. We start by extending the
definition of $L_v(s, \pi)$ to the ramified and archimedean
places. We observe that in equation (\ref{eq:gcd}), the right hand
side is indeed the greatest common denominator of the family of
rational functions $\{\zZ(W, s)\}$. Since we have already noted
that such a g.c.d. exists, even when the given representation is
not unramified, we set
\begin{equation}\label{eq:1}
L_v(s, \pi) = {\rm g.c.d.}\,\, \{\zZ(W, s)\},
\end{equation}
when $v < \infty$. Also, when $v = \infty$, we can choose an
appropriate product of Tate's archimedean L-functions, denoted by
$L_\infty(s, \pi, r)$, such that the ratio
\begin{equation}\label{eq:2}
\frac{\zZ(W, s)}{L_\infty(s, \pi)}
\end{equation}
is an entire function for all $W \in \mathcal{W}(\pi_v, \psi_v)$,
and it is a nowhere vanishing function for some choice of $W$.

With this extension, we now proceed to outline the proof. Let $S$
be a set of places, including the place at infinity, such that for
$v \notin S$, all the data is unramified. We set
\begin{equation}
L_S(s, \pi) = \prod_{v \notin S} L_v (s, \pi).
\end{equation}
For each ``ramified" non-archimedean place $p$, we choose $W_p$
such that $\zZ(W_p, s)=1$. Also for the archimedean $v$, we choose
$W_v$ such that $\zZ(W_v, s)$ is a non-vanishing entire function
$e^{g(s)}$. If we set $W = \prod_v W_v$, with $W_p = W_p^0$ for $p
\in S$, we have
\begin{equation}
\zZ(\varphi, s) = e^{g(s)}L_S(s, \pi),
\end{equation}
implying the holomorphicity of $L_S$. This immediately implies the
continuation of $L$ to a meromorphic function with only a finite
number of poles. It is this point with which the present note is
concerned.

We finally turn to the functional equation of the completed
L-function. Choosing $W_p$ so that $\zZ(W_p, s) = L_p(s, \pi)$, we
have
\begin{align*}
L(s, \pi) & = \zZ(\varphi, s) \\
& = \tilde{\zZ}(\varphi^w, 1-s) \\
& = \bigl( \prod_{p \in S} \tilde{\zZ}(W_p^w, 1-s) \bigr) L_S(1-s,
\tilde{\pi}) \\
& = \bigl( \prod_{p \in S} \frac{\tilde{\zZ}(W_p^w, 1-s)}{L_p(1-s,
\tilde{\pi})} \bigr) L(1-s, \tilde{\pi}) \\
& =\bigl( \prod_{p \in S} \frac{\gamma(\pi_p, \psi_p, s)\zZ(W_p,
s)}{L_p(1-s,\tilde{\pi})} \bigr) L(1-s, \tilde{\pi}) \\
& =\bigl( \prod_{p \in S} \frac{\gamma(\pi_p, \psi_p, s)L_p(s,
\pi)}
{L_p(1-s,\tilde{\pi})} \bigr) L(1-s, \tilde{\pi}) \\
& =\bigl( \prod_{p \in S} \epsilon(s, \pi_p, \psi_p) \bigr) L(1-s,
\tilde{\pi}),
\end{align*}
where
\begin{equation*}
\epsilon(s, \pi_p, \psi_p) = \frac{\gamma(\pi_p, \psi_p, s)L_p(s,
\pi)} {L_p(1-s,\tilde{\pi})}.
\end{equation*}
Hence, if we set
\begin{equation*}
\epsilon(s, \pi) =\prod_{p \in S} \epsilon(s, \pi_p, \psi_p),
\end{equation*}
we have the functional equation
\begin{equation}
L(s, \pi) = \epsilon(s, \pi) L(1-s, \tilde{\pi}),
\end{equation}
as anticipated by Langlands. One last note is that the function
$\epsilon(s, \pi)$ is a monomial function of $s$. In particular,
it has no poles or zeroes.
\begin{rem} In equation (\ref{eq:asymptotic}), if $c_i=\mu_i.
v^{r_i}$, with $\mu_i$ a quasi-character, we have
\begin{equation*}
L_p(s, \pi) = \prod_{i=1}^N L(s, \mu_i)^{r_i}.
\end{equation*}
The L-functions appearing on the right hand side are Tate's local
L-factors for the quasi-characters $\mu_i$. This implies that in
order to give an explicit calculations of the local L-factors, we
need to determine the finite functions $c_i$. In
\cite{Jacquet-Langlands}, this is established by a case by case
analysis of representation types for $\pi_p$, i.e. principal
series vs. special representations vs. supercuspidals.
\end{rem}

\section{Preliminaries on $\GSp(4)$}
\subsection{Bessel functionals}
We recall the notion of Bessel model introduced by Novodvorsky and
Piatetski-Shapiro \cite{Novodvorsky-PS}. We follow the exposition
of \cite{Furusawa}. Let $S \in M_2(\Q)$ be such that $S = \,^t S$.
We define the discriminant $d = d(S)$ of $S$ by $d(S) = -4 \det
S$. Let us define a subgroup $T = T_S$ of $\GL(2)$ by
\begin{equation*}
T = \{ g \in \GL(2)\, \vert \, \,^t g S g = \det g. S \}.
\end{equation*}
Then we consider $T$ as a subgroup of $\GSp(4)$ via
\begin{equation*}
t \mapsto \begin{pmatrix} t \\ & \det t. \,^t t ^{-1}
\end{pmatrix},
\end{equation*}
$t \in T$.

Let us denote by $U$ the subgroup of $\GSp(4)$ defined by
\begin{equation*}
U = \{ u(X) = \begin{pmatrix} I_2 & X \\ & I_2 \end{pmatrix} \,
\vert \, X = \,^t X \}.
\end{equation*}
Finally, we define a subgroup $R$ of $\GSp(4)$ by $R = TU$.

Let $\psi$ be a non-trivial character of $\Q \backslash \A$. Then
we define a character $\psi_S$ on $U(\A)$ by $\psi_S(u(X)) =
\psi(\tr (SX))$ for $X = \,^t X \in \M_2(\A)$. Usually when there
is no danger of confusion, we abbreviate $\psi_S$ to $\psi$. Let
$\Lambda$ be a character of $T(\Q) \backslash T(\A)$. Denote by
$\Lambda \otimes \psi_S$ the character of $R(\A)$ defined by
$(\Lambda \otimes \psi)(tu) = \Lambda(t) \psi_S(u)$ for $t \in
T(\A)$ and $u \in U(\A)$.

Let $\pi$ be an automorphic cuspidal representation of
$\GSp_4(\A)$ and $V_\pi$ its space of automorphic functions. We
assume that
\begin{equation}\label{compatible}
\Lambda \vert_{\A^\times} = \omega_\pi.
\end{equation}
Then for $\varphi \in V_\pi$, we define a function $B_\varphi$ on
$\GSp_4(\A)$ by
\begin{equation}\label{Bessel}
B_\varphi(g) = \int_{Z_\A R_\Q \backslash R_\A} (\Lambda \otimes
\psi_S)(r)^{-1}. \varphi(rh) \, dh.
\end{equation}
We say that $\pi$ has a global Bessel model of type $(S, \Lambda,
\psi)$ for $\pi$ if for some $\varphi \in V_\pi$, the function
$B_\varphi$ is non-zero. In this case, the $\C$-vector space of
functions on $\GSp_4(\A)$ spanned by $\{ B_\varphi \, \vert \,
\varphi \in V_\pi \}$ is called the space of the global Bessel
model of $\pi$.

Similarly, one can consider local Bessel models. Fix a local field
$\Q_v$. Define the algebraic groups $T_S$, $U$, and $R$ as above.
Also, consider the characters $\Lambda$, $\psi$, $\psi_S$, and
$\Lambda \otimes \psi_S$ of the corresponding local groups. Let
$(\pi, V_\pi)$ be an irreducible admissible representation of the
group $\GSp(4)$ over $\Q_v$. Then we say that the representation
$\pi$ has a local Bessel model of type $(S, \Lambda, \psi)$ if
there is a functional $\lambda_B$ on $(V_\pi^\infty)'$, a
continuous linear functional on $V_\pi^\infty$ in such a way that
\begin{equation*}
\lambda_B(\pi(r) v ) = (\Lambda \otimes \psi_S)(r)\lambda_B(v),
\end{equation*}
for all $r \in R(\Q_v)$, $v \in V_\pi$. Also, we require that
$\lambda_B$ would have some continuity properties similar to the
ones satisfied by local Whittaker functionals.

In this work, we will be interested in two different types of
Bessel models corresponding to two choices of the symmetric matrix
$S$. The two choices of $S$ are:
\begin{enumerate}
\item $S = \begin{pmatrix} & 1 \\ 1 \end{pmatrix}$,
\item $S = \begin{pmatrix} 1 \\ & d \end{pmatrix}$, with $d$ a positive
square-free rational number.
\end{enumerate}
Below, we will determine the subgroups $T_S$, and $R$, and
explicitly write down the corresponding global Bessel functionals.
We fix an irreducible automorphic cuspidal representation $\pi$ of
$\GSp_4(\A)$ and a unitary character $\psi$ of $\A$ throughout.

(1) $S = \begin{pmatrix} & 1 \\ 1
\end{pmatrix}$. This is the case of interest for us in this work. In this case, the subgroup $T_S$
is equal to the subgroup consisting of diagonal matrices. A
straightforward analysis then shows that for every character
$\Lambda$ of $T_S(\Q)\backslash T_S(\A)$ subject to
(\ref{compatible}), there is a Hecke character of $\A^\times$ such
that the global Bessel functional (\ref{Bessel}) is given by
\begin{equation*}
B^{\rm{split}}_\chi(g ; \varphi)= \int_{F^\times \backslash
\A^\times} \varphi^U
(\begin{pmatrix} y \\ & 1 \\ && 1 \\
&&& y \end{pmatrix}) \chi(y)\, d^\times y.
\end{equation*}
Here when $\phi$ is a cusp form on $\GSp(4)$, we have set
\begin{equation*}
\phi^U(g) = \int_{ (F \backslash \A)^3} \phi ( \begin{pmatrix} 1
&& u & w \\ & 1& w & v \\ &&1 \\ &&&1 \end{pmatrix}g)
\psi^{-1}(w)\, du \, dv \, dw.
\end{equation*}
(2) $S = \begin{pmatrix} 1 \\ & d \end{pmatrix}$. In this case,
the subgroup $T_S$ is equal to a non-split torus. Then there is a
Hecke character of the torus $T_S$, say $\chi$, in such a way that
\begin{equation*}
B_\chi(g ; \varphi)= \int_{T_S(F)\A^\times \backslash T_S(\A)}
\varphi^U (\begin{pmatrix} \alpha \\ \det \alpha. \,^t \alpha^{-1}
\end{pmatrix}) \chi(\alpha)\, d\alpha,
\end{equation*}
with $\phi^U$ defined as before. The case of immediate interest is
the case where $d=1$, in which case,
\begin{equation*}\begin{split}
T_S & = \{ g \in \GL_2 \, \vert \, \,^t g . g = \det g \, \} \\
& = \{ \begin{pmatrix} a & b \\ -b & a \end{pmatrix} \, \vert \,
a^2 + b^2 \in \GL_1 \}.
\end{split}\end{equation*}
The problems of existence of Bessel functionals for this choice of
the matrix $S$ seem to be more delicate. We have considered these
problems in \cite{Takloo-Bighash2}.

\subsection{Theta}\label{Section:Integral}
In this section we collect various results on theta correspondence
that we will use in the sequel. In fact, this paper is a review of
\cite{Roberts}. We have adapted the results of that paper to the
case of our interest, split orthogonal spaces of signature
$(2,2)$. Other references of interest are \cite{Harris-Kudla,
Harris-Soudry-Taylor}.

Let $V$ be the vector space $\M_2$, of the two by two matrices,
equipped with the quadratic form $\det$. Let $(,)$ be the
associated non-degenerate inner product, and $H = \GO(V, (,))$ be
the group of orthogonal similitudes of $V$, $(,)$. The group
$\GL(2) \times \GL(2)$ has a natural involution $t$ defined by
$t(g_1, g_2) = (^t b_2 ^{-1}, ^t b_1^{-1})$, where the superscript
$t$ stands for the transposition. Let $\tilde{H} = (\GL(2) \times
\GL(2)) \rtimes <t>$ be the semi-direct product of $\GL(2) \times
\GL(2)$ with the group of order two generated by $t$. There is a
sequence
\begin{equation}
1 \longrightarrow \mathbb{G}_m \longrightarrow \tilde{H}
\longrightarrow H \longrightarrow 1,
\end{equation}
where the homomorphism $\rho: \tilde{H} \rightarrow H$ is defined
by $\rho(g_1, g_2)(v) = g_1 v g_2^{-1}$, and $\rho(t) v = \,^t v$,
for all $g_1, g_2 \in \GL(2)$ and $v \in V$. Also, $\mathbb{G}_m
\rightarrow \tilde{H}$ is the natural map $z \mapsto (z,z)\times 1
$. It follows that the image of the subgroup $\GL(2) \times \GL(2)
\subset \tilde{H}$ under $\rho$ is the connected component of the
identity of $H$.

Let $F$ be a local field of characteristic zero, with $F=\R$ if
$F$ is archimedean. Fix a non-trivial unitary character $\psi$ of
$F$. The Weil representation $\omega$ of $\Sp(4, F) \times \rm{O}
(V, F)$ defined with respect to $\psi$ is the unitary
representation on $L^2(V^2)$ given by
\begin{align*}
\omega(1, h) \varphi(x) & = \varphi(h^{-1} x), \\
\omega\left( \begin{pmatrix} a \\ & \,^t a^{-1}
\end{pmatrix}\right)\varphi(x) & = \left| \det a \right|^2 \varphi(x
a), \\
\omega\left( \begin{pmatrix} 1 & b \\ & 1\end{pmatrix}
\right)\varphi(x) & =
\psi\left( \frac{1}{2} \tr (bx, x) \right)\varphi(x) , \\
\omega\left( \begin{pmatrix} & 1 \\ -1 \end{pmatrix} \right)
\varphi(x) & = \gamma \hat{\varphi}(x).
\end{align*}
Here, $\hat{\varphi}$ is the Fourier transform defined by
$$
\varphi(x) = \int_{V^2} \varphi(x') \psi(\tr (x, x')) \, dx'
$$
with $dx'$ self-dual, and $\gamma$ is a certain fourth root of
unity on $\psi$. If $h \in {\rm O}(V, F)$, $a \in \GL(2, F)$, $b
\in \M_n(F)$ with $\,^t b = b$ and $x = (x_1, x_2), x' = (x_1',
x_2') \in V^2$, we write $h^{-1} x = (h^{-1} x_1, h^{-1} x_2)$, $x
a = (x_1, x_2) ( a_{ij})$, $(x, x') = ((x_i, x_j'))$, $bx = b \,^t
(x_1, x_2)$.

If $F$ is non-archimedean, $\omega$ preserves the space
$\mathcal{S}(V^2)$; by $\omega$ we mean $\omega$ acting on the
latter space. When $F = \R$, we will work with Harish-Chandra
modules of real reductive groups. Fix $K_1 = \Sp(4, \R) \cap {\rm
O}(4, \R)$ as a maximal compact subgroup of $\Sp(4, \R)$. We
denote the Lie algebra of $\Sp(4, \R)$ by $\mathfrak{g_1}=
{\mathfrak sp}(4, \R)$. Let $V^+$ and $V^-$ be positive and
negative definite subspaces of $X$, respectively, such that $V =
V^+ \bot V^-$. Then a maximal compact subgroup of ${\rm O}(V, \R)$
is ${\rm O}(V^+, \R) \times {\rm V}(V^-, \R) \simeq {\rm O}(2,\R)
\times {\rm O}(2,\R)$. The Lie algebra of ${\rm O}(V, \R)$ is
$\mathfrak{h_1}= \mathfrak{o}(V, \R)$. Let $\mathcal{S}(V^2) =
\mathcal{S}_\psi(V^2)$ be the subspace of $L^2 (V^2)$ consisting
of the functions
$$
p(x) \exp \left[ - \frac{1}{2} |c| \left(\tr (x^+, x^+) - \tr(x^-,
x^-) \right)\right].
$$
Here $p$ is a polynomial, and $(x^+, x^+)$ and $(x^-, x^-)$ are $2
\times 2$ matrices with $(i, j)$-th entries $(x_i^+, x_j^+)$ and
$(x_i^+, x_j^+)$ respectively, where $x_i = x_i^+ + x_i^-$
corresponding to the decomposition of $V$; $c \in \R^\times$ is
such that $\psi(t) = \exp(ict)$. Then $\mathcal{S}(V^2)$ is a
$(\mathfrak{g}_1 \times \mathfrak{h}_1, K_1, J_1)$ module under
$\omega$; this is the Harish-Chandra module we will work with
throughout. Often, for the sake of uniformity in presentation, one
uses the notation and terminology of genuine representations for
archimedean places as well. The reader has to keep on mind,
however, that this is just a matter of convenience.

Let $\mathcal{R}(\rm{O}(V, F))$ be the set of elements of ${\rm
Irr}\,(\rm{O}(V, F))$ which are non-zero quotients of $\omega$,
and define $\mathcal{R}(\Sp(4, F))$ similarly. Again, the reader
will have to keep in mind that at the archimedean place, we are
working with underlying Harish-Chandra modules. Suppose $F$ is
real or non-archimedean of odd residual characteristic. Then the
set
$$
\left\{ (\pi, \sigma) \in \mathcal{R}(\Sp(4, F)) \times
\mathcal{R}(\rm{O}(V, F))\, \vert \, \Hom_{\Sp(4, F) \times
\rm{O}(V, F)} (\omega, \pi \otimes \sigma ) \ne 0 \right\}
$$
is the graph of a bijection, denoted by $\theta$ in either
direction, between the corresponding sets. When $F$ is
non-archimedean of even residual characteristic, one can establish
the same for tempered representations. We refer the reader to
\cite{Roberts}, section 1, for more information.

We now recall the extended Weil representation for similitude
groups. Define
$$
R_V(F) = \left\{ (g, h) \in \GSp(4, F) \times \GO(V, F) \, \vert
\, \nu(g) = \nu(h) \right\}.
$$
The Weil representation of $\Sp(4, F) \times \rm{O}(V, F)$ on
$L^2(V^2)$ extends to a unitary representation of $R_V(F)$ via
$$
\omega(g, h) \varphi = \left|\nu(h)\right|^{-2} \omega(g
\begin{pmatrix} 1 \\ & \nu(g) \end{pmatrix}^{-1}, 1) (\varphi \circ
h^{-1}).
$$
We would still like to consider the action of $R_V(F)$ on
$\mathcal{S}(V^2)$, but one has to take some care when considering
the archimedean place, as in this case $\mathcal{S}(V^2)$ is
preserved only at the level of Harish-Chandra modules; we refer
the reader to \cite{Roberts} for details. We denote the resulting
genuine representation of $R_V$, in the non-archimedean case, or
the $(\mathfrak{r}_\infty, L_\infty)$ Harish-Chandra module, in
the archimedean case, again by $\omega$.

In analogy with the isometry case, one can ask when
$\Hom_{R_V}(\omega, \pi \otimes \sigma) \ne 0$ for $\pi \in
\rm{Irr}\, (\GSp(4,F))$ and $\sigma \in \rm{Irr}(\rm{GO}(V, F))$.
Here $\mathcal{R}$ for each group is the collection of
representations of the similitude group which when restricted to
the corresponding isometry group have a non-zero component in
$\mathcal{R}$. Then by theorem 1.8 of \cite{Roberts}, parts 1, 3,
5, $\Hom_{R_V}(\omega, \pi \otimes \sigma) \ne 0$ defines a
bijection between $\mathcal{R}(\GSp(4, F))$ and
$\mathcal{R}(\GO(V, F))$. Again, over a non-archimedean field of
even residual characteristic one has to restrict to an appropriate
class of representations.  Again, one denotes the resulting
bijection by $\theta$. Proposition 1.11 of \cite{Roberts} states
that $\theta$ maps unramified representations to unramified
representations.

Let $(\pi_1, \pi_2)$ be a pair of representations of $\GL_2$ over
the local field $F$ with $\omega_{\pi_1}.\omega_{\pi_2} = 1$.
Roberts \cite{Roberts} has associated to $(\pi_1, \pi_2)$ an
$L$-packet in $\GSp(4)$. Essentially, the idea is to consider the
representation $\pi = \pi_1 \otimes \pi_2$ of $\GSO(V, F)$ and
then consider all possible extensions of $\pi$ to $\GO(V, F)$;
then consider the theta lifts of all such extended representations
to $\GSp(4, F)$. We describe the $L$-parameter giving this packet
in the archimedean situation. If $g_i = \begin{pmatrix} \alpha_i &
\beta_i \\ \gamma_i & \delta_i \end{pmatrix}$, $i=1,2$, we set
$$
S(g_1, g_2) = \begin{pmatrix} \alpha_1 && \beta_1& \\
                              &        \alpha_2&&\beta_2\\
                              \gamma_1&&\delta_1& \\
                              &\gamma_2&&\delta_2 \end{pmatrix}.
$$
For $i=1,2$, let $\rho_i: W_\R \to \GL_2(\C)$ be the $L$ parameter
of $\pi$. Then define an $L$-parameter $\varphi(\rho_1,\rho_2) :
W_\R \to \GSp(4, \C)$ by
$$
\varphi(\rho_1, \rho_2)(z) = S(\rho_1(z), \rho_2(z)^{-1}),
$$
$z \in W_\R$. We take for granted the fact that the $L$ packet
defined by Roberts in the archimedean situation is the $L$ packet
associated to $\varphi(\rho_1, \rho_2)$ by Langlands. We refer the
reader to section 4 of \cite{Roberts}, in particular pages 283-285
for basic properties of the $L$ packets.

We now turn our attention to global theta correspondence for the
similitude groups \cite{Roberts}, section 5. In order to define
global theta correspondence we need a global Weil representation.
Fix a non-trivial unitary character of $\A$ trivial on $\Q$. For a
place $v$ of $\Q$, let $\omega_v$ be the representation defined
above. Let $x_1, \dots, x_4$ be a vector space basis of $\M_2(\Q)$
over $\Q$. Let $(g, h) \in R_V(\A)$. Then for almost all places
$v$, $\omega_v(g_v, h_v)$ fixes the characteristic function of
$\O_v x_1 + \dots + \O_v x_4$. Let $\mathcal{S}(V(\A)^2)$ be the
restricted algebraic direct product $\otimes_v
\mathcal{S}(V(\Q_v)^2)$ which is naturally an $R_V(\A_f) \times
(\mathfrak{r}_\infty, L_\infty)$-module. For $\varphi \in
\mathcal{S}(V(\A)^2)$ and $(g, h) \in R_V(\A)$, define
$$
\theta(g, h; \varphi) = \sum_{x \in V(\Q)^2} \omega(g, h)
\varphi(x).
$$
This series converges absolutely and is left $R(\Q)$ invariant.
Fix a right invariant quotient measure on ${\rm O}(V, \Q)
\backslash {\rm O}(V, \A)$. Let $f$ be a cusp form on $\GO(V,
\A)$. For $g \in \GSp(4, \A)$ define
$$
\theta(f, \varphi)(g) = \int_{{\rm O}(V, \Q) \backslash {\rm O}(V,
\A)} \theta(g, h_1 h; \varphi)f(h_1 h)\, dh_1,
$$
where $h \in \GO(V, \A)$ is any element such that $(g, h) \in
R_V(\A)$. This integral converges absolutely, does depend on the
choice of $h$, and the function $\theta(f, \varphi)$ on $\GSp(4,
\A)$ is left $\GSp(4, \Q)$ invariant. The function $\theta(f,
\varphi)$ is an automorphic function on $\GSp(4, \A)$ of central
character equal to the central character of $f$. If $V$ is a
$\GO(V, \A) \times (h_\infty, J_\infty)$ subspace of the space of
cusp forms on $\GO(V, \A)$ of central character $\chi$, then we
denote by $\Theta(V)$ the $\GSp(4, \A_f) \times (g_\infty,
K_\infty)$ subspace of the space of automorphic forms on $\GSp(4,
\A)$ of central character $\chi$ generated by all the $\theta(f,
\varphi)$ for $f \in V$ and $\varphi \in \mathcal{S}(V(\A)^2)$.

For computational purposes, we need to make the above
considerations explicit. Here the notation may be slightly
different from above. Suppose $\pi_1$ and $\pi_2$ are two
irreducible cuspidal automorphic representations of $\GL_2(\A)$
satisfying
\begin{equation*}
\omega_{\pi_1}. \omega_{\pi_2} = 1.
\end{equation*}
Then for $\varphi_1$ and $\varphi_2$ cusp forms in the spaces of
$\pi_1$ and $\pi_2$, respectively, one can think of
\begin{equation*}
\varphi (h_1, h_2) = \varphi_1(h_1) \varphi_2(h_2),
\end{equation*}
as a cusp form on the algebraic group $\rho(\tilde{H})$. We extend
the definition of $\varphi$ to $H$ by defining it to be right
invariant under the compact totally disconnected group $<t>(\A) =
\prod_v <t>$.

Define the subgroup $H_1$ consisting of elements $(h_1, h_2)$
satisfying
\begin{equation*}
\det(h_1) = \det(h_2).
\end{equation*}
Then if $\pi_1$ and $\pi_2$ are two automorphic cuspidal
representations of the group $\GL(2)$ with
\begin{equation*}
\omega_{\pi_1}.\omega_{\pi_2}=1,
\end{equation*}
and
\begin{equation*}
\pi_1 \neq \tilde\pi_2,
\end{equation*}
then one can naturally think of the pair $(\pi_1, \pi_2)$ as an
automorphic cuspidal representation of the group $H$. If
$\varphi_1$ and $\varphi_2$ are cusp forms on $\GL_2(\A)$,
belonging to the spaces of the representations $\pi_1$ and
$\pi_2$, respectively, we define a cuspidal function
$\theta(\varphi_1, \varphi_2; f)$ on $\GSp(4, \A)$ by
\begin{equation*}
\theta(\varphi_1, \varphi_2; f)(g) = \int_{H_1(F) \backslash
H_1(\A)} \theta(g; h_1 h^1, h_2 h^2; f) \varphi_1(h_1 h^1)
\varphi_2(h_2 h^2) \, d (h_1, h_2),
\end{equation*}
where the pair $(h^1, h^2)$ is chosen such that
\begin{equation*}
\det h^1 (\det h^2)^{-1} = \nu (g).
\end{equation*}
Here $f$ is a Bruhat-Schwartz function on $\M_2(\A) \times
\M_2(\A)$, and
\begin{equation*}
\theta (g; h_1 h^1, h_2 h^2; f) = \sum_{ M_1, M_2 \in \M_2 (F)}
\omega(g; h_1h^1, h_2 h^2)f(M_1, M_2),
\end{equation*}
where $\omega$ is the Weil representation of \cite{Harris-Kudla}.
We note this is different from the definition given earlier. Let
$\Theta(\pi_1, \pi_2)$ be the vector space generated by the
functions $\theta(\varphi_1, \varphi_2; f)$ for all choices of
$\varphi_1$, $\varphi_2$, and $f$ as above. Then $\Theta(\pi_1,
\pi_2)$ is an irreducible generic automorphic cuspidal
representation of $\GSp(4)$. In fact, this is the generic element
of the global $L$ packet defined by Roberts \cite{Roberts}. If
$\Theta(\pi_1, \pi_2) = \otimes_v \Theta_v(\pi_1, \pi_2)$, then
$\Theta_v(\pi_1, \pi_2)$ depends only on the $v$ components of
$\pi_1, \pi_2$, and is the generic element of corresponding local
$L$ packet.

\subsection{The Spinor L-function for \GSp(4)} In this section, we
review the integral representation given by Novodvorsky
\cite{Novodvorsky} for $G = \GSp(4)$. The details of the material
in the following paragraphs appear in \cite{Bump},
\cite{Takloo-Bighash}.

Let $\varphi$ be a cusp form on $\GSp(4, \A)$, belonging to the
space of an irreducible cuspidal automorphic representation $\pi$.
Consider the integral
\begin{align*}
Z_N(s, \phi, \mu)= \int_{{\mathbb A}^\times /
\Q^\times}\int_{({\mathbb A}/\Q)^3} \phi\Biggl( &
\begin{pmatrix}1&x_2&x_4& \\ & 1 \\ &&1& \\ &z&-x_2&1 \end{pmatrix}
\begin{pmatrix}y \\ &y \\ &&1 \\ &&&1 \end{pmatrix} \Biggr)\\
& \times \psi(-x_2)\mu(y) \vert y \vert^{s - {1 \over 2}}\, dz \,
dx_2 \, dx_4 \, d^\times y.
\end{align*}
Since $\phi$ is left invariant under the matrix
\begin{equation*}
\begin{pmatrix} &&&1 \\ &&1& \\ &-1 \\ -1 \end{pmatrix},
\end{equation*}
this integral has a functional equation $s \rightarrow 1-s$. A
usual unfolding process as sketched in \cite {Bump} then shows
that
\begin{equation}\label{eq:fundamental}
\Z_N(s, \phi, \mu) = \int_{{\mathbb A}^\times} \int_{\mathbb A}
W_\phi
\begin{pmatrix} y \\ &y \\ &&1 \\ &x&&1 \end{pmatrix} \mu(y)
\vert y \vert ^{s - {3 \over 2}} \, dx \, d^\times y.
\end{equation}

Here the Whittaker function $W_\varphi$ is given by
\begin{align*}
W_\phi(g) = \int_{({\mathbb A} / \Q)^4} \phi \Biggl(
\begin{pmatrix} 1 & x_2 \\ &1 \\ && 1 \\ &&-x_2&1 \end{pmatrix} &
\begin{pmatrix} 1 &&x_4 &x_3 \\ &1&x_3&x_1 \\ &&1 \\ &&&1 \end{pmatrix}
g \Biggr) \\
& \times \psi^{-1}(x_1 + x_2) \, dx_1 \, dx_2 \, dx_3 \, dx_4
\end{align*}
Equation (\ref{eq:fundamental}) implies that, in order for
$Z_N(\varphi, s)$ to be non-zero, we need to assume that
$W_\varphi$ is not identically equal to zero. A representation
satisfying this condition is called ``generic." Every irreducible
cuspidal representation of $\GL(2)$ is generic. On other groups,
however, there may exist non-generic cuspidal representations. In
fact, those representations of $\GSp(4)$ which correspond to
holomorphic cuspidal Siegel modular forms are not generic.

{\em From this point on, we assume that all the representations of
\GSp(4), local or global, which appear in the text are generic.}

If $\varphi$ is chosen correctly, the Whittaker function may be
assumed to decompose locally as $W(g) = \prod_v W_v(g_v)$, a
product of local Whittaker functions. Hence, for $\Re s$ large, we
obtain
\begin{equation}
\zZ(\varphi, s) = \prod_v \zZ(W_v, s),
\end{equation}
where
\begin{equation}
Z_N(W_v, s) =\int_{F_v^\times} \int_{F_v} W_v \Bigl(
\begin{pmatrix} y  \\ & y \\ & & 1 \\ &x&&1 \end{pmatrix}
\Bigr) \vert y \vert^{s - {3 \over 2}}\, dx \, d^\times y.
\end{equation}
As usual, we have a functional equation: There exists a
meromorphic function $\gamma(\pi_v, \psi_v, s)$ (rational function
in $\N v^{-s}$ when $v < \infty$) such that
\begin{equation}
Z_N(W_v, s) = \gamma(\pi_v, \psi_v, s)\tilde{\zZ}(W_v^w, 1-s),
\end{equation}
with $w$ as above,
\begin{equation*}
\tilde{\zZ}(W_v, s) =\int_{F_v^\times} \int_{F_v} W_v \Bigl(
\begin{pmatrix} y  \\ & y \\ & & 1 \\ &x&&1 \end{pmatrix}
\Bigr) \chi_v^{-1}(y) \vert y \vert^{s - {3 \over 2}}\, dx \,
d^\times y,
\end{equation*}
and $\chi_v$ the central character of $\pi_v$.

We also consider the unramified calculations. Suppose $v$ is any
nonarchimedean place of $F$ such that $W_v$ is right invariant by
$\GSp(4, \O_v)$ and such that the largest fractional ideal on
which $\psi_v$ is trivial is $\O$. Then the Casselman-Shalika
formula \cite{Casselman-Shalika} allows us to calculate the last
integral (cf. \cite{Bump}). The result is the following:
\begin{equation}
\zZ(W_v, s) = L(s, \pi_v, {\rm Spin}).
\end{equation}
Let us explain the notation. The connected L-group $^L G^0$ is
$\GSp(\C)$. Let $^L T$ be the maximal torus of elements of the
form
\begin{equation*}
t(\alpha_1, \alpha_2, \alpha_3, \alpha_4) = \begin{pmatrix}
\alpha_1 \\ & \alpha_2 \\ && \alpha_3 \\ &&& \alpha_4
\end{pmatrix},
\end{equation*}
where $\alpha_1 \alpha_4 = \alpha_2 \alpha_3$. The fundamental
dominant weights of the torus are $\lambda_1$ and $\lambda_2$
where
\begin{equation*}
\lambda_1 t(\alpha_1, \alpha_2, \alpha_3, \alpha_4) = \alpha_1,
\end{equation*}
and
\begin{equation*}
\lambda_2 t(\alpha_1, \alpha_2, \alpha_3, \alpha_4) = \alpha_1
\alpha_3^{-1}.
\end{equation*}
The dimensions of the representation spaces associated with these
dominant weights are four and five, respectively. In our notation,
Spin is the representation of $\GSp(4, \C)$ associated with the
dominant weight $\lambda_1$, i.e. the standard representation of
$\GSp(4, \C)$ on $\C^4$. The L-function $L(s, \pi, {\rm Spin})$ is
called the Spinor, or simply the Spin, L-function of $\GSp(4)$.

Next step is to use the integral introduced above to extend the
definition of the Spinor L-function to ramified non-archimedean
and archimedean places.

We now sketch the computation of the local non-archimedean Euler
factors of the Spin L-function given by the integral
representation introduced above. In order for this to make sense,
we need the following lemma:
\begin{lem}[Theorem 2.1 of \cite{Takloo-Bighash}]\label{lem:bound}
Suppose $\Pi$ is a generic representation of $\GSp(4)$ over a
non-archimedean local field $K$, $q$ order of the residue field.
For each $W \in \W(\Pi, \psi)$, the function $\zZ(W, s)$ is a
rational function of $q^{-s}$, and the ideal $\{\zZ(W, s) \}$ is
principal.
\end{lem}
\noindent {\em Sketch of proof.} For $W \in \W(\Pi, \psi)$, we set
\begin{equation}
Z(W, s) = \int_K W \left( \begin{pmatrix} y \\ & y \\ && 1 \\ &&&
1
\end{pmatrix} \right) \vert y \vert^{s - {3 \over 2}}\, d^\times y.
\end{equation}
The first step of the proof is to show that the vector space $\{
Z(W, s) \}$ is the same as $\{\zZ(W, s) \}$ (cf. Proposition 3.2
of \cite{Takloo-Bighash}). Next, we use the asymptotic expansions
of the Whittaker functions along the torus to prove the existence
of the g.c.d. for the ideal $\{ Z(W, s) \}$. Indeed, Proposition
3.5 of \cite{Takloo-Bighash} (originally a theorem in
\cite{Casselman-Shalika}) states that there is a finite set of
finite functions $S_\Pi$, depending only on $\Pi$, with the
following property: for any $W \in \W(\Pi, \psi)$, and $c \in
S_\Pi$, there is a Schwartz-Bruhat function $\Phi_{c, W}$ on $K$
such that
\begin{equation*}
W \Bigl( \begin{pmatrix} y \\ & y \\ && 1 \\ &&& 1 \end{pmatrix}
\Bigr) = \sum_{c \in S_\Pi} \Phi_{c, W}(y) c(y) \vert y \vert^{3
\over 2}.
\end{equation*}
The lemma is now immediate. $\square$

We have the following theorem:
\begin{thm}\label{thm:main} Suppose $\Pi$ is a generic representation of the
group $\GSp(4)$ over a non-archimedean local field $K$. Then
\begin{enumerate}
\item If $\Pi$ is supercuspidal, or is a sub-quotient of a
representation induced from a supercuspidal representation of the
Klingen parabolic subgroup, then $L(s, \pi, {\rm Spin})=1$. \item
If $\pi$ is a supercuspidal representation of $\GL(2)$ and $\chi$
a quasi-character of $K^\times$, and $\Pi = \pi \rtimes \chi$ is
irreducible, we have
\begin{equation*}
L(s, \Pi, {\rm Spin}) = L(s, \chi).L(s, \chi.\omega_\pi).
\end{equation*}
\item If $\chi_1$, $\chi_2$, and $\chi_3$ are quasi-characters of
$K^\times$, and $\Pi = \chi_1 \times \chi_2 \rtimes \chi_3$ is
irreducible, we have
\begin{equation*}
L(s, \Pi, {\rm Spin}) = L(s, \chi_3).L(s, \chi_1 \chi_3).L(s,
\chi_2 \chi_3).L(s, \chi_1 \chi_2 \chi_3).
\end{equation*}
\item \label{classification} When $\Pi$ is not irreducible, one
can prove similar statements for the generic subquotients of $\Pi
= \pi \rtimes \chi$ (resp. $\Pi = \chi_1 \times \chi_2 \rtimes
\chi_3$) according to the classification theorems of Sally-Tadic
\cite{Sally-Tadic} and Shahidi \cite{Shahidi2} (cf. theorems 4.1
and 5.1 of \cite{Takloo-Bighash}).
\end{enumerate}
\end{thm}
\begin{rem} Sally and Tadic \cite{Sally-Tadic} and Shahidi
\cite{Shahidi2} have completed the classification of
representations supported in the Borel and Siegel parabolic
subgroups. In particular, they have determined for which
representations the parabolic induction is reducible. From their
result, one can immediately establish a classification for all the
generic representations supported in the Borel or Siegel parabolic
subgroups.
\end{rem}
\noindent {\em Sketch of proof.} By the proof of the lemma, we
need to determine the asymptotic expansion of the Whittaker
functions in each case. The argument consists of several steps:
\begin{step} {\em Bound the size of $S_\Pi$.} Fix $c \in S_\Pi$,
and define a functional $\Lambda_c$ on $\W(\Pi, \psi)$ by
\begin{equation}
\Lambda_c(W) = \Phi_{c, W}(0).
\end{equation}
If $c, c' \in S_\Pi$, and $c \ne c'$, the two functionals
$\Lambda_c$ and $\Lambda_{c'}$ are linearly independent.
Furthermore, the functionals $\Lambda_c$ belong to the dual of a
certain twisted Jacquet module $\Pi_{N, \bar{\theta}}$ (notation
from \cite{Takloo-Bighash}, page 1095). Hence $\# S_\Pi = \dim
\Pi_{N, \bar{\theta}}$. Then one uses an argument similar to those
of \cite{Shalika}, distribution theory on p-adic manifolds, to
bound the dimension of the Jacquet module. The result (proposition
3.9 of \cite{Takloo-Bighash}) is that if $\Pi$ is supercuspidal or
supported in the Klingen parabolic subgroup (resp. Siegel
parabolic, resp. Borel parabolic), then $ \# S_\Pi = 0$ (resp.
$\leq 2$, resp. $\leq 4$). Note that this already implies the
first part of the theorem.
\end{step}
From this point on, we concentrate on the Siegel parabolic
subgroup, the Borel subgroup case being similar. We fix some
notation. Suppose $\Pi = \pi \rtimes \chi$, with $\pi$
supercuspidal of $\GL(2)$. Let $\lambda_\Pi$ (resp. $\lambda_\pi$)
be the Whittaker functional of $\Pi$ (resp. $\pi$) from
\cite{Shahidi}. It follows from the proof of the lemma
\ref{lem:bound} that, for $f \in \Pi$, there is a positive number
$\delta(f)$, such that
\begin{equation*}
\lambda_\Pi( \Pi \Bigl( \begin{pmatrix} y \\ & y \\ && 1 \\ &&& 1
\end{pmatrix} \Bigr) f)= \sum_{c \in S_\Pi} \Lambda_c(f) c(y)
\vert y \vert^{3 \over 2},
\end{equation*}
for $\vert y \vert < \delta(f)$. Here, $\Lambda_c$ is the obvious
functional on the space of $\Pi$.
\begin{step}{\em Uniformity.} For $f \in Ind ( \pi \times \chi
\vert P \cap K, K)$, and $\tau \in \C$, define $f_\tau$ on $G$ by
\begin{equation*}
f_\tau(pk) = \delta_P(p)^{\tau + {1 \over 2}} \pi \otimes \chi (p)
f(k).
\end{equation*}
It is clear that $f_\tau$ is a well-defined function on $G$, and
that it belongs to the space of a certain induced representation
$\Pi_\tau$. The {\em Uniformity Theorem} (Proposition 3.9 of
\cite{Takloo-Bighash}) asserts that one can take $\delta(f_\tau) =
\delta(f)$.
\end{step}
\begin{step} {\em Regular representations.} This is
the case where $\omega_\pi \ne 1$. In this situation, we have
\begin{equation}\begin{split}
\lambda_\Pi& ( \Pi \Bigl(  \begin{pmatrix} y \\ & y \\ && 1 \\ &&&
1
\end{pmatrix} \Bigr)  f) = \\   & \lambda_\pi (A(w, \Pi)(f)(e)) \chi(y)
\vert y \vert^{3 \over 2} + C( w\Pi, w^{-1})^{-1}
\lambda_\pi(f(e)) \chi(y)\omega_\pi(y)\vert y \vert^{3 \over 2},
\end{split}\end{equation}
for $\vert y \vert < \delta(f)$. Here $w = \begin{pmatrix} && 1
\\ &&&1 \\ -1 \\ & -1 \end{pmatrix}$, $A(w, \Pi)$ is the
intertwining integral of \cite{Shahidi}, and $C(w\Pi, w^{-1})$ is
the local coefficient of \cite{Shahidi}. The proof of this
identity follows from the the above lemma \ref{lem:bound}, and the
{\em Multiplicity One Theorem} \cite{Shalika}. The idea is to find
one term of the asymptotic expansion using the open cell; then
apply the long intertwining operator to find the other term.
\end{step}
Note that the identity of {\em Step} 3 also applies to reducible
cases. For example, if $f \in \Pi$ is in the kernel of the
intertwining operator $A(w, \Pi)$, the first term of the right
hand side vanishes.
\begin{step} {\em Irregular Representations.} The idea is the
following: we twist everything in {\em Step} 3 by the complex
number $\tau$, so that the resulting representation $\Pi_\tau$ is
regular. By {\em Step} 2, the identity still holds uniformly for
all $\tau$. By a theorem of Shahidi \cite{Shahidi} (essentially
due to Casselman and Shalika \cite{Casselman-Shalika}), we know
that the left hand side of the identity is an entire function of
$\tau$. This implies that the poles of the right hand side, coming
from the intertwining operator and the local coefficient, must
cancel out. Next, we let $\tau \to 0$. An easy argument
(l'Hopital's rule!) shows the appearance of $\chi(y) \vert y
\vert^{3 \over 2}$ and $\chi(y) \vert y \vert^{3 \over 2} \log_q
\vert y \vert$ in the asymptotic expansion.
\end{step}
This finishes the sketch of proof of the theorem. $\square$
\begin{coro}\label{very-highly-ramified}
Let $\pi$ be an irreducible generic representation of $\GSp(4)$
over a non-archimedean local field $K$. Let $\mu$ be a
quasi-character of $K^\times$. If $\mu$ is highly ramified, we
have
$$
L(s, \pi \otimes \mu) = 1.
$$
\end{coro}
\section{Bessel Functionals and Integral Representations}
In the global situation, there is a simple relationship between
the integral representation of the previous section and split
Bessel functionals. The following simple observation which for the
ease of reference we separate as a lemma forms the fundamental
idea of this paper:
\begin{lem}\label{Bessel-Zeta} We have
\begin{equation*}\begin{split}
B^{\rm{split}}_{\mu \vertdot^{s-\frac{1}{2}}}(I_4;
\phi)=\int_{{\mathbb A}^\times} \int_{\mathbb A} W_\phi\left(
\begin{pmatrix} y \\ &y \\ &&1 \\ &x&&1 \end{pmatrix}w^{-1}\right) \mu(y)
\vert y \vert ^{s - {3 \over 2}} \, dx \, d^\times y,
\end{split}\end{equation*}
with
\begin{equation*}
w = \begin{pmatrix} 1 \\ &&&1 \\ &&1 \\ &-1 \end{pmatrix}.
\end{equation*}
\end{lem}
This motivates the following definition.

\begin{defn} For $\varphi_1$, $\varphi_2$, and $f$ as above and $\mu$ a
Hecke character, we define
\begin{align*}
\zZ(\varphi_1, \varphi_2, f; \mu) &=
B^{\rm{split}}_{\mu\vertdot^{-\frac{1}{2}}}(I_4; \theta(\varphi_1,
\varphi_2; f)) \\
& = \int_{F^\times \backslash \A^\times} \theta(\varphi_1,
\varphi_2; f)^U( \begin{pmatrix} y \\ & 1 \\ && 1 \\
&&& y \end{pmatrix}) \mu(y) \vert y \vert^{- {1 \over 2}}\,
d^\times y.
\end{align*}
\end{defn}
Here if $\phi$ is a cusp form on $\GSp(4)$, we have set
\begin{equation*}
\phi^U(g) = \int_{ (F \backslash \A)^3} \phi ( \begin{pmatrix} 1
&& u & w \\ & 1& w & v \\ &&1 \\ &&&1 \end{pmatrix}g)
\psi^{-1}(w)\, du \, dv \, dw.
\end{equation*}
We prove that the above integral is an infinite product of local
integrals. We do so by finding an expression relating our function
$\zZ(\varphi_1, \varphi_2, f; s)$ to the Jacquet-Langlands zeta
functions of $\varphi_1$, and $\varphi_2$.

Before stating our proposition, we recall a notation from
\cite{Jacquet-Langlands}. If $\phi$ is a cusp form on
$\GL_2(\A_F)$, in the space of a representation $\pi$, $\mu$ a
Hecke character, and $h \in \GL_2(\A_F)$, we set
\begin{equation*}
Z(\phi, h, \mu) = \int_{F^\times \backslash \A^\times} \phi(
\begin{pmatrix} a \\ & 1 \end{pmatrix} h ) \mu(a) \vert a
\vert^{-\frac{1}{2}} \, d^\times a,
\end{equation*}
and
\begin{equation*}
\tilde{Z}(\phi, h, \mu) = \int_{F^\times \backslash \A^\times}
\phi(\begin{pmatrix} a \\ & 1 \end{pmatrix} h )
\omega_\pi(a)^{-1}\mu(a) \vert a \vert^{-\frac{1}{2}} \, d^\times
a
\end{equation*}
Then, we have the following proposition:
\begin{prop} For $\varphi_1$, $\varphi_2$, and $f$ as above,
we have
\begin{equation*}\begin{split}
\zZ(\varphi_1, \varphi_2, f; \mu) = \int_{D(\A) \backslash
H_1(\A)}
& Z(\varphi_1, h_1, \mu) Z(\varphi_2, h_2, \mu^{-1} \vertdot ) \\
& L(h_1 , h_2) f (\begin{pmatrix} 1 & 0
\\ 0 & 0 \end{pmatrix},
\begin{pmatrix} 0 & 0 \\ 0 & 1 \end{pmatrix} )
dh_1 \, dh_2
\end{split}\end{equation*}
\end{prop}
\begin{proof} First, we obtain an expression for $\theta(\varphi_1,
\varphi_2; f)^U$. We start by the following:
\begin{align*}
\theta  (\varphi_1, \varphi_2 & ; f)^U(g) \\
 = & \int_{ (F \backslash \A)^3} \theta(\varphi_1, \varphi_2; f)
( \begin{pmatrix} 1 && u & w \\ & 1& w & v \\ &&1 \\
&&&1 \end{pmatrix}g) \psi^{-1}(w)\, du \, dv \, dw
\end{align*}
\begin{align*}
 = & \int_{ (F \backslash \A)^3} \int_{H_1(F) \backslash H_1(\A)}
\theta( \begin{pmatrix} 1 && u & w \\ & 1& w & v \\ &&1 \\ &&&1
\end{pmatrix} g; h_1 h^1, h_2 h^2; f) \\
& \varphi_1(h_1 h^1) \varphi_2(h_2 h^2) \, d (h_1, h_2)
\psi^{-1}(w)\, du \, dv \, dw,
\end{align*}
where $h^1$ and $h^2$ are chosen in such a way that
\begin{equation*}
\det h^1. (\det h^2)^{-1} = \nu (g).
\end{equation*}
Next, it follows from the definition of $\theta$ that
\begin{equation}\begin{split}\label{thetaU}
\theta  (\varphi_1& , \varphi_2 ; f)^U(g) = \\
& \int_{H_1(F) \backslash H_1(\A)}\varphi_1( h_1 h^1)
\varphi_2(h_2 h^2) G_f(h_1 h^1, h_2 h^2; g)\, dh_1 \, dh_2,
\end{split}\end{equation}
where
\begin{align*}
G_f & (h_1h^1, h_2 h^2 ; g) = \\
& \sum_{M_1, M_2} \int_{(F \backslash
\A)^3} \omega(\begin{pmatrix} 1 && u & w \\ & 1& w & v \\ &&1 \\
&&&1 \end{pmatrix} g, h_1 h^1, h_2 h^2) f(M_1, M_2) \\
&  \psi^{-1}(w) \, du \, dv \, dw.
\end{align*}
Next, for fixed $M_1$ and $M_2$ we have
\begin{align*}
&\int_{(F \backslash \A)^3} \omega(\begin{pmatrix} 1 && u & w \\ & 1& w & v \\ &&1 \\
&&&1 \end{pmatrix} g, h_1 h^1, h_2 h^2) f(M_1, M_2) \psi^{-1}(w)
\, du \, dv \, dw \\
& = \omega( g , h_1 h^1, h_2 h^2) f(M_1, M_2)\\
&  \int_{(F \backslash \A)^3} \psi( tr \begin{pmatrix} u & w \\ w
& v \end{pmatrix}
\begin{pmatrix} \det M_1 & B(M_1, M_2) - \frac{1}{2} \\
B(M_2, M_1) - \frac{1}{2} & \det M_2 \end{pmatrix}) \, du \, dv \,
dw.
\end{align*}
Next, we have the following straightforward lemma:
\begin{lem} For any $2 \times 2$ matrix $A \in \M_2(\A)$, we have
\begin{equation*}
\int_{(F \backslash \A)^3} \psi (tr \begin{pmatrix} u & w \\ w & v
\end{pmatrix} A)\, du \, dv \, dw =0,
\end{equation*}
unless $A = \begin{pmatrix} 0 & 0 \\ 0 & 0 \end{pmatrix}$, in
which case the value of the integral is equal to 1.
\end{lem}
The lemma implies that
\begin{equation*}
G_f(h_1 h^1, h_2 h^2; g) = \sum_{(M_1, M_2) \in \mathcal{S}}
\omega(g, h_1 h^1, h_2 h^2) f(M_1, M_2),
\end{equation*}
where
\begin{equation*}
\mathcal{S} = \{ (X, Y) \in \M_2(F) \times \M_2(F) \, \vert \,
\det X =0, \det Y = 0, \det (X+Y)=1 \}.
\end{equation*}
\begin{lem} The set $\mathcal{S}$ consists of a single orbit under
the action of $H_1(F)$. The point $P= ( \begin{pmatrix} 1 & 0
\\ 0 & 0 \end{pmatrix}, \begin{pmatrix} 0 & 0 \\ 0 & 1
\end{pmatrix} )$ belongs to $\mathcal{S}$. The stabilizer of $P$
in $H_1(F)$ is the subgroup $D(F)$.
\end{lem}
Consequently,
\begin{equation*}\begin{split}
G_f(h_1 h^1, & h_2 h^2; g) =  \\
& \sum_{ \gamma \in D(F) \backslash H_1(F)} \omega(1, \gamma)
\omega (g, h_1 h^1, h_2 h^2) f (
\begin{pmatrix} 1 & 0
\\ 0 & 0 \end{pmatrix}, \begin{pmatrix} 0 & 0 \\ 0 & 1
\end{pmatrix} ).
\end{split}\end{equation*}
Inserting the right hand side of this expression for $G_f$ in
equation (\ref{thetaU}) gives
\begin{equation}\begin{split}\label{thetaU2}
\theta & (\varphi_1 , \varphi_2 ; f)^U(g) = \\
& \int_{D(F) \backslash H_1(\A)}\varphi_1( h_1 h^1) \varphi_2(h_2
h^2) \omega (g, h_1 h^1, h_2 h^2) f (\begin{pmatrix} 1 & 0
\\ 0 & 0 \end{pmatrix}, \begin{pmatrix} 0 & 0 \\ 0 & 1
\end{pmatrix} ) \, dh_1 \, dh_2,
\end{split}\end{equation}
We now turn our attention to $\zZ(\varphi_1, \varphi_2, f; s)$.
For this purpose, we need to first simplify $\omega (g, h_1 h^1,
h_2 h^2) f (\begin{pmatrix} 1 & 0
\\ 0 & 0 \end{pmatrix}, \begin{pmatrix} 0 & 0 \\ 0 & 1
\end{pmatrix} )$, when $g = {\tiny \begin{pmatrix} y \\ & 1 \\ && 1
\\&&&y \end{pmatrix}}$, $h^1 = \begin{pmatrix} y \\ & 1
\end{pmatrix}$, and $h^2 = \rm{ identity}$, say. We have
\begin{equation*}
\omega(\begin{pmatrix} y \\ & 1 \\ && 1
\\&&&y \end{pmatrix}, h_1 \begin{pmatrix} y \\ & 1
\end{pmatrix}, h_2) f (\begin{pmatrix} 1 & 0
\\ 0 & 0 \end{pmatrix}, \begin{pmatrix} 0 & 0 \\ 0 & 1
\end{pmatrix} )
\end{equation*}
\begin{equation*}
= \omega(\begin{pmatrix} y \\ & 1 \\ && 1
\\&&&y \end{pmatrix}\begin{pmatrix} 1 \\ & 1 \\ && y^{-1}
\\&&&y^{-1} \end{pmatrix})L(h_1 \begin{pmatrix} y \\ & 1
\end{pmatrix}, h_2) f (\begin{pmatrix} 1 & 0
\\ 0 & 0 \end{pmatrix}, \begin{pmatrix} 0 & 0 \\ 0 & 1
\end{pmatrix} )
\end{equation*}
\begin{equation*}
= \vert y \vert^2 L(h_1 \begin{pmatrix} y \\ & 1
\end{pmatrix}, h_2) f (\begin{pmatrix} y & 0
\\ 0 & 0 \end{pmatrix}, \begin{pmatrix} 0 & 0 \\ 0 & 1
\end{pmatrix} )
\end{equation*}
\begin{equation*}
= f (\begin{pmatrix} y^{-1} & 0 \\ 0 & 1 \end{pmatrix} h_1^{-1}
\begin{pmatrix} y & 0 \\ 0 & 0 \end{pmatrix} h_2,
\begin{pmatrix} y^{-1} & 0 \\ 0 & 1 \end{pmatrix} h_1^{-1}
\begin{pmatrix} 0 & 0 \\ 0 & 1 \end{pmatrix} h_2)
\end{equation*}
\begin{equation*}
= f (\begin{pmatrix} y^{-1} & 0 \\ 0 & 1 \end{pmatrix} h_1^{-1}
\begin{pmatrix} y & 0 \\ 0 & 1 \end{pmatrix}
\begin{pmatrix} 1 & 0 \\ 0 & 0 \end{pmatrix} h_2,
\begin{pmatrix} y^{-1} & 0 \\ 0 & 1 \end{pmatrix} h_1^{-1}
\begin{pmatrix} y & 0 \\ 0 & 1 \end{pmatrix}
\begin{pmatrix} 0 & 0 \\ 0 & 1 \end{pmatrix} h_2).
\end{equation*}
Hence, for the choices of $g$, $h^1$, and $h^2$ as above, we have
\begin{align*}
\omega (g, h_1 h^1, h_2 h^2) & f (\begin{pmatrix} 1 & 0
\\ 0 & 0 \end{pmatrix}, \begin{pmatrix} 0 & 0 \\ 0 & 1
\end{pmatrix} ) = \\
&  L(\begin{pmatrix} y^{-1} & 0 \\ 0 & 1 \end{pmatrix} h_1
\begin{pmatrix} y & 0 \\ 0 & 1 \end{pmatrix}, h_2) f (\begin{pmatrix} 1 & 0
\\ 0 & 0 \end{pmatrix}, \begin{pmatrix} 0 & 0 \\ 0 & 1
\end{pmatrix} ).
\end{align*}
This equation combined with equation (\ref{thetaU2}) gives
\begin{equation*}\begin{split}
\theta & (\varphi_1 , \varphi_2 ; f)^U (\begin{pmatrix} y \\ & 1
\\ && 1 \\&&&y \end{pmatrix}) =  \int_{D(F) \backslash
H_1(\A)}\varphi_1( h_1 \begin{pmatrix} y
\\ & 1 \end{pmatrix}) \varphi_2(h_2) \\
& L(\begin{pmatrix} y^{-1} & 0 \\ 0 & 1
\end{pmatrix} h_1
\begin{pmatrix} y & 0 \\ 0 & 1 \end{pmatrix}, h_2) f (\begin{pmatrix} 1 & 0
\\ 0 & 0 \end{pmatrix}, \begin{pmatrix} 0 & 0 \\ 0 & 1
\end{pmatrix} ) \, dh_1 \, dh_2.
\end{split}\end{equation*}
Next, we make a change of variables
\begin{equation*}
(h_1, h_2) \mapsto (\begin{pmatrix} y \\ & 1 \end{pmatrix} h_1
\begin{pmatrix} y^{-1} \\ & 1 \end{pmatrix}, h_2)
\end{equation*}
to obtain
\begin{equation*}\begin{split}
\theta & (\varphi_1 , \varphi_2 ; f)^U(
\begin{pmatrix} y \\ & 1 \\ && 1 \\&&&y \end{pmatrix}) = \\
& \int_{D(F) \backslash H_1(\A)}\varphi_1( \begin{pmatrix} y
\\ & 1 \end{pmatrix} h_1) \varphi_2(h_2)
L(h_1 , h_2) f (\begin{pmatrix} 1 & 0
\\ 0 & 0 \end{pmatrix}, \begin{pmatrix} 0 & 0 \\ 0 & 1
\end{pmatrix} ) \, dh_1 \, dh_2.
\end{split}\end{equation*}
Next,
\begin{align*}
\zZ (\varphi_1& , \varphi_2, f; \mu) \\
 =& \int_{F^\times \backslash \A^\times}
\theta  (\varphi_1 , \varphi_2 ; f)^U(\begin{pmatrix} y \\ & 1 \\
&& 1 \\&&&y \end{pmatrix})\mu(y)\vert y \vert^{ - \frac{1}{2}}\,
d^\times y \\
=& \int_{F^\times \backslash \A^\times} \int_{D(F) \backslash
H_1(\A)}\varphi_1( \begin{pmatrix} y \\ & 1
\end{pmatrix} h_1) \varphi_2(h_2 ) \\
& L(h_1 , h_2) f (\begin{pmatrix} 1 & 0 \\ 0 & 0
\end{pmatrix}, \begin{pmatrix} 0 & 0 \\ 0 & 1 \end{pmatrix} )
\mu(y)\vert y \vert^{- \frac{1}{2}}\, dh_1 \, dh_2 \, d^\times y.
\end{align*}
At this stage, we use the obvious isomorphism
\begin{equation*}
F^\times \backslash \A^\times \longrightarrow D(F) \backslash
D(\A),
\end{equation*}
given by
\begin{equation*}
a \mapsto (\begin{pmatrix} a \\ & 1 \end{pmatrix},
\begin{pmatrix} a \\ & 1 \end{pmatrix})
\end{equation*}
to obtain
\begin{align*}
\zZ (\varphi_1& , \varphi_2, f; \mu) \\
=&\int_{F^\times \backslash \A^\times} \int_{D(\A) \backslash
H_1(\A)} \int_{F^\times \backslash \A^\times} \varphi_1(
\begin{pmatrix} y \\ & 1 \end{pmatrix}
\begin{pmatrix} a \\ & 1 \end{pmatrix}
h_1) \varphi_2(
\begin{pmatrix} a \\ & 1 \end{pmatrix}h_2 ) \\
& L(\begin{pmatrix} a \\ & 1 \end{pmatrix}h_1 ,
\begin{pmatrix} y \\ & 1 \end{pmatrix}h_2) f
(\begin{pmatrix} 1 & 0 \\ 0 & 0
\end{pmatrix},
\begin{pmatrix} 0 & 0 \\ 0 & 1 \end{pmatrix} )
\mu(y)\vert y \vert^{ - \frac{1}{2}}\, d^\times a \, dh_1 \, dh_2
\,
d^\times y \\
=&\int_{F^\times \backslash \A^\times} \int_{D(\A) \backslash
H_1(\A)} \int_{F^\times \backslash \A^\times} \varphi_1(
\begin{pmatrix} ya \\ & 1 \end{pmatrix}
h_1) \varphi_2(
\begin{pmatrix} a \\ & 1 \end{pmatrix}h_2 ) \\
& L(h_1 , h_2) f (\begin{pmatrix} 1 & 0
\\ 0 & 0
\end{pmatrix},
\begin{pmatrix} 0 & 0 \\ 0 & 1 \end{pmatrix} )
\mu(y)\vert y \vert^{ - \frac{1}{2}}\, d^\times a \, dh_1 \, dh_2
\,
d^\times y \\
=&\int_{F^\times \backslash \A^\times} \int_{D(\A) \backslash
H_1(\A)} \int_{F^\times \backslash \A^\times} \varphi_1(
\begin{pmatrix} y \\ & 1 \end{pmatrix}
h_1) \varphi_2(
\begin{pmatrix} a \\ & 1 \end{pmatrix}h_2 ) \\
& L(h_1 , h_2) f (\begin{pmatrix} 1 & 0
\\ 0 & 0 \end{pmatrix},
\begin{pmatrix} 0 & 0 \\ 0 & 1 \end{pmatrix} )
\mu(y)\vert y \vert^{ - \frac{1}{2}} \mu^{-1}(a)\vert a
\vert^{\frac{1}{2}}\, d^\times a \, dh_1 \, dh_2 \, d^\times y,
\end{align*}
after a change of variable $y \mapsto ya^{-1}$. The proposition
now follows from a simple re-arrangement of the last expression.
\end{proof}
\subsection{The zeta integral of two complex variables; Euler product}
In order to study the zeta integral $\zZ(\varphi_1, \varphi_2, f;
\mu)$, we would have liked to introduce a function of two complex
variables $s_1$, $s_2$ as follows: For $\varphi_1$, $\varphi_2$,
and $f$ as above, and $\mu$ Hecke character, we set
\begin{equation*}\begin{split}
\zZ(\varphi_1, \varphi_2, f; \mu, \vertdot^{s_1}, \vertdot^{s_2})
= \int_{D(\A) \backslash
H_1(\A)} & Z(\varphi_1, h_1,\mu \vertdot^{s_1} ) Z(\varphi_2, h_2,\mu^{-1}\vertdot^{s_2} ) \\
& L(h_1 , h_2) f (\begin{pmatrix} 1 & 0
\\ 0 & 0 \end{pmatrix},
\begin{pmatrix} 0 & 0 \\ 0 & 1 \end{pmatrix} )
dh_1 \, dh_2,
\end{split}\end{equation*}
with $s_1, s_2 \in \C$. Unfortunately, however, this integral is
not well-defined for $s_2 \ne 1 - s_1$. In order to circumvent
this problem we proceed as follows.

If $\phi$ is a cusp form on $\GL_2(\A_F)$, we define its Whittaker
function by
\begin{equation*}
W_\phi(g) = \int_{F \backslash \A} \phi( \begin{pmatrix} 1 & x \\
& 1 \end{pmatrix} g) \psi(x)^{-1}\, dx,
\end{equation*}
for $g \in \GL_2(\A_F)$. Then, we have the Fourier expansion
\begin{equation*}
\phi(g) = \sum_{\alpha \in F^\times} W_\phi(\begin{pmatrix} \alpha
\\ & 1 \end{pmatrix} g),
\end{equation*}
with the right hand side a uniformly convergent series on compact
sets in $\GL_2(A)$. It is then a classical observation of
\cite{Jacquet-Langlands} that for $\Re s$ large, we have
\begin{equation*}
Z(\phi, h, \mu\vertdot^s) = \int_{\A} W_\phi( \begin{pmatrix} a \\
& 1
\end{pmatrix} h ) \mu(a)\vert a \vert^{s - \frac{1}{2}}\, d^\times a.
\end{equation*}
We denote the right hand side of this equation by $Z(W_\phi, h ,
s)$.

We have a formal identity as follows:
\begin{equation*}\begin{split}
\zZ(\varphi_1, \varphi_2, f; \mu, \vertdot^{s_1}, \vertdot^{s_2})
= \int_{D(\A) \backslash H_1(\A)} & Z(W_{\varphi_1}, h_1,\mu
\vertdot^{s_1} ) Z(W_{\varphi_2},
h_2,\mu^{-1}\vertdot^{s_2} ) \\
& L(h_1 , h_2) f (\begin{pmatrix} 1 & 0
\\ 0 & 0 \end{pmatrix},
\begin{pmatrix} 0 & 0 \\ 0 & 1 \end{pmatrix} )
dh_1 \, dh_2.
\end{split}\end{equation*}

Next, we consider the Euler product. We choose $\varphi_i$, for
$i=1,2$, so that
\begin{equation*}
W_{\varphi_i} = \otimes_{v \in \mathcal{M}_F} W_v^i.
\end{equation*}
Also, we choose $f$ to be a pure tensor of the form
\begin{equation*}
\otimes_{v \in \mathcal{M}_F} f_v,
\end{equation*}
with $f_v$ unramified for almost all $v$.

With this choice of the data, we have yet another formal identity
\begin{equation}
\zZ(\varphi_1, \varphi_2, f; \mu, \vertdot^{s_1}, \vertdot^{s_2})
=\prod_{v \in \mathcal{M}_F} \zZ_v( W_v^1, W_v^2, f_v;\mu_v,
\vertdot_v^{s_1}, \vertdot_v^{s_2} ).
\end{equation}
Here, we have set
\begin{equation*}\begin{split}
\zZ_v(W_v^1, W_v^2, f_v; \mu_v, \vertdot^{s_1}, \vertdot^{s_2}) =
\int_{D(F_v) \backslash
H_1(F_v)} & Z(W_v^1, h_1, \mu_v \vertdot_v^{s_1}) Z(W_v^2, h_2, \mu_v^{-1}\vertdot_v^{s_2}) \\
& L(h_1 , h_2) f_v (\begin{pmatrix} 1 & 0
\\ 0 & 0 \end{pmatrix},
\begin{pmatrix} 0 & 0 \\ 0 & 1 \end{pmatrix} )
dh_1 \, dh_2.
\end{split}\end{equation*}
Also, for $W_v$ a Whittaker function on a local group
$\GL_2(F_v)$, and $h \in \GL_2(F_v)$, we have used the notation
$Z(W_v, h, \mu_v)$ to denote
\begin{equation*}
\int_{F_v^\times} W_v (\begin{pmatrix} a \\ & 1 \end{pmatrix} h )
\mu_v(a) \vert a \vert^{ - \frac{1}{2}}\, d^\times a.
\end{equation*}
The idea is to make sense out of the expression for
\begin{equation*}
\zZ_v(W_v^1, W_v^2, f_v; \mu_v, \vertdot^{s_1}, \vertdot^{s_2})
\end{equation*}
for $\Re s_1, \Re s_2$ large. For this we use the following lemma:
\begin{lem}\label{measure}
Let $v \in \mathcal{M}_F$, and $\Psi$ a continuous function of
compact support on $D(F_v) \backslash H_1(F_v)$. Choose an
arbitrary lift $\Phi'$ of $\Phi$ to $\GL_2(F_v) \times
\GL_2(F_v)$. The functional $\mu(\Phi)$ defined by
\begin{equation*}
\int_{K_v} \int_{F_v^2} \int_{F_v^\times} \Phi'(
\begin{pmatrix} 1 & u \\ & 1 \end{pmatrix} k_1,
\begin{pmatrix} \epsilon \\ & \epsilon^{-1} \end{pmatrix}
\begin{pmatrix} 1 & v \\ & 1 \end{pmatrix} k_2)|\epsilon|^{-1}\,
d^\times \epsilon \, du \, dv \, dk_1 \, dk_2,
\end{equation*}
for an appropriate choice of a local maximal compact (and open for
$v$ non-archimedean), defines an invariant measure on $D(F_v)
\backslash H_1(F_v)$. Furthermore, this measure has the following
property: Fix a Haar measure $\mu_D$ on $D(F_v)$, and for any
continuous function of compact support $\Psi$ on $H_1(F_v)$, set
\begin{equation*}
\Psi_D(x) = \int_{D(F_v)} \Psi( y x ) \, d \mu_1(y),
\end{equation*}
for $x \in D(F_v) \backslash H_1(F_v)$. Then the functional
$\mu_2$ defined by
\begin{equation*}
\mu_2(\Psi) = \mu( \Psi_D),
\end{equation*}
with $\Psi$ as above defines a Haar measure on $H_1(F_v)$.
\end{lem}
\begin{defn}
We set
\begin{equation*}\begin{split}
\zZ_v& (W_v^1, W_v^2, f_v; \mu_v, \vertdot^{s_1}, \vertdot^{s_2})
\\=& \int_{u, v \in F_v}\int_{\epsilon \in F_v^\times}\int_{K_v^2}
f( k_1^{-1} \begin{pmatrix} \epsilon^{-1} & \epsilon^{-1} v \\ 0 &
0
\end{pmatrix} k_2, k_1^{-1} \begin{pmatrix} 0 & - u \epsilon \\ 0
& \epsilon \end{pmatrix} k_2) \\
& \omega_{\pi_2}(\epsilon) \vert \epsilon \vert^{2 s_2 -2} \bigl(
\int_{F_v^\times} W_1 ( \begin{pmatrix} \alpha \\ & 1
\end{pmatrix} k_1) \e( u \alpha) \mu(\alpha) \vert \alpha \vert ^{s_1 - {1
\over 2}}\, d
^\times \alpha \bigr) \\
&\bigl( \int_{F_v^\times} W_2 ( \begin{pmatrix} \beta \\ & 1
\end{pmatrix} k_2) \e( v \beta) \mu^{-1}(\beta)\vert \beta \vert ^{s_2 - {1
\over 2}}\, d ^\times \beta \bigr)\, du \, dv \, d^\times \epsilon
\, dk_1 \, dk_2.
\end{split}\end{equation*}
\end{defn}
We immediately observe that if the integral is convergent, it is
well-defined.
\begin{prop} Suppose $W_1, W_2$ are two Whittaker functions of
$\GL_2(F_v)$ belonging to the spaces of representations $\pi_1,
\pi_2$, respectively, with $\omega_{\pi_1}. \omega_{\pi_2} =1$.
Then the integral $\zZ(W_1, W_2, f; \mu_v, \vertdot_v^{s_1},
\vertdot_v^{s_2})$ converges absolutely for $\Re s_1, \Re s_2 \gg
0$.
\end{prop}
\begin{proof} We give a complete proof only for the case where $v$
is a real place, the proof of the non-archimedean statement being
identical. Also it is clear that we may assume that the
quasi-character $\mu_v$ is trivial. By definition, we need to show
that the integral
\begin{equation*}\begin{split}
& \int_{u, v \in \R}\int_{\epsilon \in \R_+^\times}\int_{K_v^2} f(
k_1^{-1} \begin{pmatrix} \epsilon^{-1} & \epsilon^{-1} v \\ 0 & 0
\end{pmatrix} k_2, k_1^{-1} \begin{pmatrix} 0 & - u \epsilon \\ 0
& \epsilon \end{pmatrix} k_2) \\
& \omega_{\pi_2}(\epsilon) \vert \epsilon \vert^{2 s_2 -2} \bigl(
\int_{\R^\times} W_1 ( \begin{pmatrix} \alpha \\ & 1 \end{pmatrix}
k_1) \e( u \alpha) \vert \alpha \vert ^{s_1 - {1 \over 2}}\, d
^\times \alpha \bigr) \\
&\bigl( \int_{\R^\times} W_2 ( \begin{pmatrix} \beta \\ & 1
\end{pmatrix} k_2) \e( v \beta) \vert \beta \vert ^{s_2 - {1
\over 2}}\, d ^\times \beta \bigr)\, du \, dv \, d^\times \epsilon
\, dk_1 \, dk_2.
\end{split}\end{equation*}
converges absolutely. By lemma 8.3.3 of \cite{JPSS}, there are
gauge functions $\xi_1, \xi_2$ such that
\begin{equation*}
\vert W_1 \vert \leq \xi_1, \text{ and } \vert W_2 \vert \leq
\xi_2.
\end{equation*}
This implies that
\begin{equation*}
\int_{\R^\times} \vert W_1 ( \begin{pmatrix} \alpha \\
& 1 \end{pmatrix} k_1) \e( u \alpha) \vert \alpha \vert ^{s_1 - {1
\over 2}}\vert \, d ^\times \alpha \leq \int_{\R^\times} \xi_1 (
\begin{pmatrix} \alpha \\ & 1 \end{pmatrix} )
\vert \alpha \vert ^{\sigma_1 - {1 \over 2}}\, d ^\times \alpha,
\end{equation*}
and
\begin{equation*}
\int_{\R^\times} \vert W_2 ( \begin{pmatrix} \beta \\
& 1 \end{pmatrix} k_2) \e( v \beta) \vert \beta \vert ^{s_2 - {1
\over 2}}\vert \, d ^\times \beta \leq \int_{\R^\times} \xi_2 (
\begin{pmatrix} \beta \\ & 1 \end{pmatrix} )
\vert \beta \vert ^{\sigma_2 - {1 \over 2}}\, d ^\times \beta.
\end{equation*}
The latter integrals converge absolutely for $\sigma_1, \sigma_2$
large. In order to conclude the proof, we need to study the
convergence of
\begin{equation*}\begin{split}
& \int_{u, v \in \R}\int_{\epsilon \in \R_+^\times}\int_{K_v^2} f(
k_1^{-1} \begin{pmatrix} \epsilon^{-1} & \epsilon^{-1} v \\ 0 & 0
\end{pmatrix} k_2, k_1^{-1} \begin{pmatrix} 0 & - u \epsilon \\ 0
& \epsilon \end{pmatrix} k_2) \\
& \omega_{\pi_2}(\epsilon) \vert \epsilon \vert^{2 s_2 -2}\, du \,
dv \, d^\times \epsilon \, dk_1 \, dk_2.
\end{split}\end{equation*}
We claim that this integral converges absolutely for all values of
$s_2$. In fact, if $f \in \mathcal{S}(\M_2(\R) \times \M_2(\R))$,
the function $g$ defined by
\begin{equation*}
g(X, Y) = \int_{K_v^2} f (k_1^{-1} X k_2, k_1^{-1} Y k_2)\, dk_1
\, dk_2
\end{equation*}
is in $\mathcal{S}(\M_2(\R) \times \M_2(\R))$. Thus, we must show
that
\begin{equation*}\begin{split}
& \int_{u, v \in \R}\int_{\epsilon \in \R_+^\times} f(
\begin{pmatrix} \epsilon^{-1} & \epsilon^{-1} v \\ 0 & 0
\end{pmatrix}, \begin{pmatrix} 0 & - u \epsilon \\ 0
& \epsilon \end{pmatrix} )  \omega_{\pi_2}(\epsilon) \vert
\epsilon \vert^{2 s_2 -2}\, du \, dv \, d^\times \epsilon
\end{split}\end{equation*}
converges absolutely for all $s_2$. The first observation, due to
Weil, is that the absolute value of a Schwartz-Bruhat function is
bounded by a Schwartz-Bruhat function. Consequently, we can assume
that $f$ is a positive Schwartz-Bruhat function. But now it is
clear that the function $\Xi$ defined by
\begin{equation*}
\Xi(\epsilon) =   \int_{u, v \in \R} f(
\begin{pmatrix} \epsilon^{-1} & \epsilon^{-1} v \\ 0 & 0
\end{pmatrix}, \begin{pmatrix} 0 & - u \epsilon \\ 0
& \epsilon \end{pmatrix} )\, du \, dv
\end{equation*}
is in the space $\mathcal{S}(\R^\times)$. Since our original
integral is bounded by
\begin{equation*}
\int_{\R} \Xi(\epsilon) \omega_{\pi_2}(\epsilon) \vert \epsilon
\vert^{2 \sigma_2 - 2} \, d^\times \epsilon,
\end{equation*}
the proposition is immediate.
\end{proof}
Then we have the following proposition:
\begin{prop}\label{separation}
Let $v$ be a non-archimedean place. Let $W_1$ and $W_2$ be given.
Then there is a choice of $f$ such that
$$
\zZ(W_1, W_2, f; \mu, \vertdot_v^{s_1}, \vertdot_v^{s_2}) = Z(W_1,
\mu \vertdot_v^{s_1})Z(W_2, \mu^{-1}\vertdot_v^{s_2}).
$$
\end{prop}
\begin{proof}
Let $M$ be a very large positive integer. Let $f = g \otimes h$ be
a Schwartz function such that
$$
{\rm Support}\, g \subset \begin{pmatrix} 1 \\ & 0 \end{pmatrix} +
\begin{pmatrix} \mathfrak{p}^M & \mathfrak{p}^M \\ \mathfrak{p}^M
& \mathfrak{p}^M \end{pmatrix},
$$
and
$$
{\rm Support}\, h \subset \begin{pmatrix} 0 \\ & 1 \end{pmatrix} +
\begin{pmatrix} \mathfrak{p}^M & \mathfrak{p}^M \\ \mathfrak{p}^M
& \mathfrak{p}^M \end{pmatrix}.
$$
Then upon setting,
$$
h_1 = \begin{pmatrix} 1 & - u \\ & 1 \end{pmatrix} \begin{pmatrix}
\alpha & \beta \\ \gamma & \delta \end{pmatrix}^{-1},
$$
$$
h_2 = \begin{pmatrix} \epsilon \\ & \epsilon^{-1} \end{pmatrix}
\begin{pmatrix} 1 & v \\ & 1 \end{pmatrix} \begin{pmatrix} a & b
\\ c & d \end{pmatrix},
$$
we get
$$
f \left( \begin{pmatrix} \alpha \epsilon ( a + vc )& \alpha
\epsilon ( b + vd) \\ \gamma \epsilon (a + vc) & \gamma \epsilon (
b+ vd) \end{pmatrix},
\begin{pmatrix}  c \epsilon^{-1} ( \alpha u + \beta)& d
\epsilon^{-1} ( \alpha u + \beta) \\ c \epsilon^{-1} ( \gamma u +
\delta) & d \epsilon^{-1} ( \gamma u + \delta)
\end{pmatrix} \right) \ne 0.
$$
With the choice of $f$, it is not hard to draw the following
conclusions:
\begin{enumerate}
\item $\gamma, c \in \mathfrak{p}^M$,

\item $u, v$ are integral,

\item $\epsilon$ is a unit,

\item $b + vd, \alpha u + \beta \in \mathfrak{p}^M$,

\item $\alpha \epsilon a, d \epsilon^{-1}\delta \in 1 +
\mathfrak{p}^M$.
\end{enumerate}
Next,
\begin{equation*}
Z(W_1, h_1, \mu_1 \vertdot_v^{s_1}) = \int_{\Q_v^\times} W_1
\left( \begin{pmatrix} x \\ & 1 \end{pmatrix}
\begin{pmatrix} 1 & -u \\ & 1 \end{pmatrix}
\begin{pmatrix} \alpha & \beta \\ \gamma & \delta \end{pmatrix}^{-1}\right)
\mu(x) \vert x \vert^{s_1 -\frac{1}{2}} \, d^\times x;
\end{equation*}
but
\begin{equation*}\begin{split}
\begin{pmatrix} 1 & -u \\ & 1 \end{pmatrix} \begin{pmatrix} \alpha & \beta \\ \gamma & \delta
\end{pmatrix}^{-1}=&
\begin{pmatrix} \alpha^{-1} \\ & \alpha (\alpha \delta - \beta \gamma)^{-1}
\end{pmatrix} \\ & \times \begin{pmatrix} 1 & -( \beta + u \alpha) \alpha (\alpha \delta - \beta \gamma)^{-1} \\
& 1  \end{pmatrix}
\begin{pmatrix} 1 \\ - \alpha ^{-1}\gamma & 1 \end{pmatrix},
\end{split}\end{equation*}
implying that for $M$ large, we have
\begin{equation*}\begin{split}
Z(W_1, h_1, \mu \vertdot_v^{s_1}) & = \int_{\Q_v^\times} W_1
\left( \begin{pmatrix} x  \alpha^{-1}
\\ & \alpha (\alpha \delta - \beta \gamma)^{-1}
\end{pmatrix}\right) \mu_1(x)
\vert x \vert^{s_1 -\frac{1}{2}} \, d^\times x \\
& = (\omega_{\pi_1} \mu)(\alpha^2 (\alpha \delta - \beta
\gamma)^{-1}) Z(W_1, \mu \vertdot_v^{s_1}).
\end{split}\end{equation*}
Similarly, for $M$ large,
\begin{equation*}
Z(W_1, h_2, \mu^{-1} \vertdot_v^{s_2}) = \mu^{-1}(\epsilon^{-1} d
(ad - bc)^{-1}) (\omega_{\pi_2}\mu^{-1})(\epsilon^{-1} d ) Z(W_2,
\mu^{-1} \vertdot_v^{s_2}).
\end{equation*}
The proposition is now immediate.
\end{proof}
\begin{coro}\label{coro:constant}
There is a choice of $W_1, W_2, f$ such that $$\zZ(W_1, W_2, f;
\mu, \vertdot_v^{s_1}, \vertdot_v^{s_2}) \equiv 1.$$
\end{coro}
When $W_1, W_2$ are spherical, the situation is particularly nice:
\begin{prop}\label{unramified-holo}
Suppose $v$ is a non-archimedean place, and $\pi_1, \pi_2$ are
spherical representations of $\GL_2(F_v)$ with $\omega_{\pi_1}.
\omega_{\pi_2}=1$. Also, suppose that $W_i \in \mathcal{W}(\pi_i,
\psi)$, $i=1, 2$, is the normalized $K_v$-fixed vector.
Furthermore, let $f$ be the characteristic function of $\M_2(\O_v)
\times \M_2(\O_v)$. Then for unramified quasi-character $\mu$ we
have
\begin{equation*}
\zZ(W_1, W_2, f; \mu, \vertdot_v^{s_1}, \vertdot_v^{s_2}) =
L_v(s_1, \pi_1, \mu)L(s_2, \pi_2, \mu^{-1}).
\end{equation*}
\end{prop}
\begin{proof} In order to see this, we need to verify that if
\begin{equation*}
L(h_1, h_2) f ( \begin{pmatrix} 1 & 0 \\ 0 & 0 \end{pmatrix},
\begin{pmatrix} 0 & 0 \\ 0 & 1 \end{pmatrix}) \ne 0,
\end{equation*}
for $(h_1, h_2) \in H_1(F_v)$, we must have $(h_1, h_2) \in
D(F_v)(\GL_2(\O_v) \times \GL_2(\O_v))$. For this, we start by the
observation that one can take as a set $\mathcal{R}$ of
representatives for
\begin{equation*}
D(F_v)\backslash H_1(F_v)/(\GL_2(\O_v) \times \GL_2(\O_v)),
\end{equation*}
the set of pairs of the form
\begin{equation*}
( \begin{pmatrix} 1 & u \\ & 1 \end{pmatrix}, \begin{pmatrix}
\epsilon \\ & \epsilon^{-1} \end{pmatrix} \begin{pmatrix} 1 & v \\
& 1 \end{pmatrix}).
\end{equation*}
Hence, we need to verify our claim only for elements $(h_1, h_2)$
of the above form. We have
\begin{equation*}
L(h_1, h_2) f ( \begin{pmatrix} 1 & 0 \\ 0 & 0 \end{pmatrix},
\begin{pmatrix} 0 & 0 \\ 0 & 1 \end{pmatrix}) =
f ( \begin{pmatrix} \epsilon & \epsilon v \\ 0 & 0 \end{pmatrix},
\begin{pmatrix} 0 & - u \epsilon^{-1} \\ 0 & \epsilon^{-1}
\end{pmatrix}).
\end{equation*}
Since $f$ is the characteristic function of $\M_2(\O_v) \times
\M_v(\O_v)$, for this last expression to be non-zero, we must have
$\epsilon^{\pm 1} \in \O_v$, $\epsilon v \in \O_v$, and
$\epsilon^{-1} u \in \O_v$. This in turn implies that $\epsilon
\in \O_v^\times$, and $u, v \in \O_v$. Now an application of lemma
\ref{measure} gives the result.
\end{proof}

We can now proceed to collect information about the analytic
properties of our two variable zeta function. we prove the
following proposition:
\begin{prop}\label{quotient-holo} For $W_1$, $W_2$ Whittaker functions,
and $f$ as above, the function $\zZ(W_1, W_2,
f;\mu,\vertdot^{s_1}, \vertdot^{s_2})$ has an analytic
continuation to a meromorphic function on $\C^2$. Furthermore, the
ratio
\begin{equation*}
\Psi(W_1, W_2, f; \mu, \vertdot_v^{s_1},  \vertdot_v^{s_2})=
\frac{\zZ(W_1, W_2, f; \mu, \vertdot_v^{s_1},
\vertdot_v^{s_2})}{L(s_1, \pi_1, \mu )L(s_2, \pi_2, \mu^{-1} )}
\end{equation*}
extends to an entire function on the entire $\C^2$. There is a
choice of $W_1$, $W_2$, and $f$ such that the above ratio is a
nowhere vanishing entire function.
\end{prop}
\begin{proof}
We prove only the analyticity statement; the non-vanishing follows
from proposition \ref{separation} and the corresponding $\GL(2)$
statement. We write out the details for the archimedean place. For
simplicity, we will assume that $\pi_1$ and $\pi_2$ are
irreducible principal series representations. Also we will assume
that the quasi-character $\mu$ is trivial. By lemma \ref{measure},
we need to consider the integral
\begin{equation}\begin{split}\label{entire}
& \int_{u, v \in \R}\int_{\epsilon \in \R_+^\times}\int_{K_v^2} f(
k_1^{-1} \begin{pmatrix} \epsilon^{-1} & \epsilon^{-1} v \\ 0 & 0
\end{pmatrix} k_2, k_1^{-1} \begin{pmatrix} 0 & - u \epsilon \\ 0
& \epsilon \end{pmatrix} k_2) \\
& \omega_{\pi_2}(\epsilon) \vert \epsilon \vert^{2 s_2 -2} \bigl(
\int_{\R^\times} W_1 ( \begin{pmatrix} \alpha \\ & 1 \end{pmatrix}
k_1) \e( u \alpha) \vert \alpha \vert ^{s_1 - {1 \over 2}}\, d
^\times \alpha \bigr) \\
&\bigl( \int_{\R^\times} W_2 ( \begin{pmatrix} \beta \\ & 1
\end{pmatrix} k_2) \e( v \beta) \vert \beta \vert ^{s_2 - {1
\over 2}}\, d ^\times \beta \bigr)\, du \, dv \, d^\times \epsilon
\, dk_1 \, dk_2.
\end{split}\end{equation}
For this purpose, we use the description of the Whittaker model of
a principal series representation from \cite{Jacquet-Langlands},
page 101-102. Suppose $\pi_1= \pi(\mu_1, \mu_2)$, and $\pi_2 =
\pi(\mu_3, \mu_4)$. Then there is a Schwartz function $P_i(x, y)$,
$i=1,2$, such that $W_1 = W_{P_i}$ by the following recipe. Let
\begin{equation*}
f_1(g) = (\mu_1 \nu^{1\over 2})(\det g) \int_{\R^\times} P_1 [ (0,
1) \gamma g ] (\mu_1 \mu_2^{-1}\nu) (\gamma) \, d^\times \gamma,
\end{equation*}
and
\begin{equation*}
f_2(g) = (\mu_3 \nu^{1\over 2})(\det g) \int_{\R^\times} P_2 [ (0,
1) \delta g ] (\mu_3 \mu_4^{-1}\nu) (\delta) \, d^\times \delta,
\end{equation*}
when the integrals converge. Next, we set for $i=1,2$
\begin{equation*}
W_{P_i}(g) = \int_\R f_{P_i}( \begin{pmatrix} & 1 \\ -1
\end{pmatrix} \begin{pmatrix} 1 & x \\ & 1 \end{pmatrix} g) \e(x)
\, dx.
\end{equation*}
In particular,
\begin{equation*}\begin{split}
W_{P_1}& ( \begin{pmatrix} \alpha \\ & 1 \end{pmatrix} k_1) = \\
&\int_\R \int_{\R^\times} (\mu_1 v^{1 \over 2})(\alpha) (\mu_1
\mu_2^{-1}\nu)(\gamma) P_1((-\alpha \gamma, - x \gamma)k_1) \e(x)
\, dx \, d^\times \gamma,
\end{split}\end{equation*}
and
\begin{equation*}\begin{split}
W_{P_2}& ( \begin{pmatrix} \beta \\ & 1 \end{pmatrix} k_2) = \\
&\int_\R \int_{\R^\times} (\mu_3 v^{1 \over 2})(\beta) (\mu_3
\mu_4^{-1}\nu)(\delta) P_2((-\beta \delta, - y \delta)k_2) \e(y)
\, dy \, d^\times \delta.
\end{split}\end{equation*}
These integrals may not converge, but they have analytic
continuations to entire functions of the characters $\mu_i$, $i=1,
\dots, 4$.

We need a lemma/notation:
\begin{lem} Suppose $P_1$, $P_2$, and $f$ are Schwartz-Bruhat
functions as above. Then the function $\Gamma$ whose value at
$$ (X, Y, m, n, p, q) \in \M_2(\R) \times \M_2(\R) \times \R^4$$
is given by
\begin{equation*}\begin{split}
\Gamma (X, & Y, m, n, p, q) = \\ & \int_{K^2} f(k_1^{-1} X k_2,
k_1^{-1} Y k_2) P_1((m, n) k_1) P_2((p, q)k_2)\, dk_1 \, dk_2
\end{split}\end{equation*}
is a Schwartz-Bruhat function.
\end{lem}
The integral (\ref{entire}) is now equal to
\begin{equation*}\begin{split}
&\int_{\alpha \in \R^\times} \int_{\beta \in \R^\times}
\int_{\gamma \in \R^\times} \int_{\delta \in \R^\times}
\int_{\epsilon \in \R_+^\times} \int_{u \in \R} \int_{v \in \R}
\int_{x \in \R} \int_{y \in \R} \\
&\Gamma (\begin{pmatrix}
\epsilon^{-1} & \epsilon^{-1}v \\ 0 & 0 \end{pmatrix},
\begin{pmatrix} 0 & - u \epsilon \\ 0 & \epsilon \end{pmatrix}, -
\alpha \gamma, - x \gamma, - \beta \delta, - y \delta) \\
& \omega_{\pi_2}(\epsilon)\vert \epsilon \vert^{2 s_2 -2} \e(u
\alpha) \vert \alpha \vert^{s_1 - {1 \over 2}} \e(v \beta) \vert
\beta \vert^{s_2 - {1 \over 2}} (\mu_1 v^{1 \over 2})(\alpha) \\
& (\mu_1 \mu_2^{-1}\nu)(\gamma) \e (x) (\mu_3 v^{1 \over
2})(\beta)(\mu_3 \mu_4^{-1}\nu)(\delta) \e(y) \\
& dy \, dx \, dv \, du \, d^\times \epsilon \, d^\times \delta \,
d^\times \gamma \, d^\times \beta \, d^\times \alpha.
\end{split}\end{equation*}
\begin{equation}\begin{split}\label{mess}
=&\int_{\alpha \in \R^\times} \int_{\beta \in \R^\times}
\int_{\gamma \in \R^\times} \int_{\delta \in \R^\times}
\int_{\epsilon \in \R_+^\times} \int_{u \in \R} \int_{v \in \R}
\int_{x \in \R} \int_{y \in \R} \\
&\Gamma (\begin{pmatrix} \epsilon^{-1} & \epsilon^{-1}v \\ 0 & 0
\end{pmatrix},
\begin{pmatrix} 0 & - u \epsilon \\ 0 & \epsilon \end{pmatrix}, -
\alpha \gamma, - x \gamma, - \beta \delta, - y \delta) \\
& \omega_{\pi_2}(\epsilon)\vert \epsilon \vert^{2 s_2 -2} \e(u
\alpha) \vert \alpha \vert^{s_1} \e(v \beta) \vert
\beta \vert^{s_2} (\mu_1)(\alpha) \\
& (\mu_1 \mu_2^{-1}\nu)(\gamma) \e (x) (\mu_3)(\beta)
(\mu_3 \mu_4^{-1}\nu)(\delta) \e(y) \\
& dy \, dx \, dv \, du \, d^\times \epsilon \, d^\times \delta \,
d^\times \gamma \, d^\times \beta \, d^\times \alpha.
\end{split}\end{equation}
We will abbreviate the inner $\Gamma$-expression appearing above
to
\begin{equation*}
\Gamma (\epsilon^{-1}, \epsilon^{-1}v , - u \epsilon , \epsilon ,
- \alpha \gamma, - x \gamma, - \beta \delta, - y \delta).
\end{equation*}
Next we consider the integral
\begin{equation*}\begin{split}
& \int_{u \in \R} \int_{v \in \R} \int_{x \in \R} \int_{y \in \R}
\Gamma (\epsilon^{-1}, \epsilon^{-1}v , - u \epsilon , \epsilon ,
- \alpha \gamma, - x \gamma, - \beta \delta, - y \delta) \\
& \e(x) \e(y) \e(u \alpha) \e(v\beta) \, dy \, dx \, dv \, du \\
= & \,\, \vert \gamma \vert^{-1} \vert \delta \vert^{-1} \int_{u
\in \R} \int_{v \in \R} \int_{x \in \R} \int_{y \in \R} \Gamma
(\epsilon^{-1}, v , u  , \epsilon ,
- \alpha \gamma, x , - \beta \delta, y ) \\
& \e(-{x \over \gamma}) \e(-{y \over \delta}) \e(-u {\alpha \over
\epsilon} ) \e(v\beta \epsilon) \, dy \, dx \, dv \, du \\
=& \,\, \vert \gamma \vert^{-1} \vert \delta \vert^{-1}
\widetilde{\Gamma} (\epsilon^{-1}, -\beta \epsilon , \alpha
\epsilon^{-1} , \epsilon , - \alpha \gamma, \gamma^{-1} , - \beta
\delta, \delta^{-1} ),
\end{split}\end{equation*}
where $\widetilde{\Gamma}$ is the appropriate Fourier transform of
$\Gamma$.

Going back to (\ref{mess}), we obtain
\begin{equation*}\begin{split}
&\int_{\alpha \in \R^\times} \int_{\beta \in \R^\times}
\int_{\gamma \in \R^\times} \int_{\delta \in \R^\times}
\int_{\epsilon \in \R_+^\times} \\
& \vert \gamma \vert^{-1} \vert \delta \vert^{-1}
\widetilde{\Gamma} (\epsilon^{-1}, -\beta \epsilon , \alpha
\epsilon^{-1} , \epsilon , - \alpha \gamma, \gamma^{-1} , - \beta
\delta, \delta^{-1} ) \\
& \omega_{\pi_2}(\epsilon)\vert \epsilon \vert^{2 s_2 -2} \vert
\alpha \vert^{s_1} \vert \beta \vert^{s_2} \mu_1(\alpha) (\mu_1
\mu_2^{-1} \nu)(\gamma) \mu_3(\beta)
(\mu_3 \mu_4^{-1}\nu)(\delta)  \\
&  d^\times \epsilon \, d^\times \delta \, d^\times \gamma \,
d^\times \beta \, d^\times \alpha.
\end{split}\end{equation*}
\begin{equation*}\begin{split}
= &\int_{\alpha \in \R^\times} \int_{\beta \in \R^\times}
\int_{\gamma \in \R^\times} \int_{\delta \in \R^\times}
\int_{\epsilon \in \R_+^\times} \\
& \widetilde{\Gamma} (\epsilon^{-1}, -\beta \epsilon , \alpha
\epsilon^{-1} , \epsilon , - \alpha \gamma^{-1}, \gamma , - \beta
\delta^{-1}, \delta ) \\
& \omega_{\pi_2}(\epsilon)\vert \epsilon \vert^{2 s_2 -2} \vert
\alpha \vert^{s_1} \vert \beta \vert^{s_2} \mu_1(\alpha) (\mu_1
\mu_2^{-1})(\gamma^{-1}) \mu_3(\beta)
(\mu_3 \mu_4^{-1})(\delta^{-1})  \\
&  d^\times \epsilon \, d^\times \delta \, d^\times \gamma \,
d^\times \beta \, d^\times \alpha.
\end{split}\end{equation*}
\begin{equation*}\begin{split}
= &\int_{\alpha \in \R^\times} \int_{\beta \in \R^\times}
\int_{\gamma \in \R^\times} \int_{\delta \in \R^\times}
\int_{\epsilon \in \R_+^\times} \\
& \widetilde{\Gamma} (\epsilon^{-1}, -\beta \delta \epsilon ,
\alpha \gamma
\epsilon^{-1} , \epsilon , - \alpha , \gamma , - \beta , \delta ) \\
& \omega_{\pi_2}(\epsilon)\vert \epsilon \vert^{2 s_2 -2} \vert
\alpha \vert^{s_1} \vert \gamma \vert^{s_1} \vert \beta
\vert^{s_2} \vert \delta \vert^{s_2}\mu_1(\alpha) \mu_2(\gamma)
\mu_3(\beta)
\mu_4(\delta)  \\
&  d^\times \epsilon \, d^\times \delta \, d^\times \gamma \,
d^\times \beta \, d^\times \alpha
\end{split}\end{equation*}
\begin{equation}\begin{split}\label{reference}
= &\int_{\alpha \in \R^\times} \int_{\beta \in \R^\times}
\int_{\gamma \in \R^\times} \int_{\delta \in \R^\times}
\int_{\epsilon \in \R_+^\times} \\
& \widetilde{\Gamma} (\epsilon^{-1}, -\beta \delta \epsilon ,
\alpha \gamma
\epsilon^{-1} , \epsilon , - \alpha , \gamma , - \beta , \delta ) \\
& (\mu_1 \nu^{s_1})(\alpha) (\mu_2 \nu^{s_1})(\gamma) (\mu_3 \nu^{s_2})(\beta)
(\mu_4 \nu^{s_2})(\delta) (\omega_{\pi_2} \nu^{2s_2 -2})(\epsilon)   \\
&  d^\times \epsilon \, d^\times \delta \, d^\times \gamma \,
d^\times \beta \, d^\times \alpha
\end{split}\end{equation}
after obvious changes of variables, and simple re-arrangement of
terms.

Our result now follows from the following standard lemma:
\begin{lem}
Let $\Phi$ be a Schwartz-Bruhat function on $\R^n$. Suppose
$\gamma_1, \dots, \gamma_n$ are quasi-characters. Define the
function $Z(s_1, \dots, s_n)= Z(\Phi; \gamma_1, \dots, \gamma_n;
s_1, \dots, s_n)$ of the complex variables $s_1, \dots, s_n$ by
\begin{equation*}
Z(s_1, \dots, s_n) = \int_{(\R^\times)^n} \Phi(\alpha_1, \dots,
\alpha_n) \prod_i \gamma_i(\alpha_i) \vert \alpha_i \vert^{s_i}\,
d^\times \alpha_i,
\end{equation*}
whenever the integral converges. Then the integral converges for
$\Re s_i$ large enough, for $i=1, \dots, n$. The ratio
\begin{equation*}
\frac{Z(\Phi; \gamma_1, \dots, \gamma_n; s_1, \dots, s_n)}{
\prod_{i=1}^n L(s_i, \gamma_i)}
\end{equation*}
extends to an entire function. If $\Phi \in \mathcal{S}( \R^\times
\times \R^{n-1})$, then the ratio
\begin{equation*}
\frac{Z(\Phi; \gamma_1, \dots, \gamma_n; s_1, \dots, s_n)}{
\prod_{i=2}^n L(s_i, \gamma_i)}
\end{equation*}
has an analytic continuation to an entire function.
\end{lem}
\end{proof}
\begin{coro}\label{highly-ramified}
Let $v$ be a non-archimedean place. Then in the above situation
for $\mu$ highly ramified $\zZ(W_1, W_2, f; \mu, \vertdot_v^{s_1},
\vertdot_v^{s_2})$ extends to an entire function of $s_1, s_2$.
\end{coro}
\begin{coro}
Let $W_1, W_2$ be flat sections of Whittaker spaces as in the last
section. Then the function $\Psi(W_1, W_2, f; \mu,
\vertdot_v^{s_1},  \vertdot_v^{s_2})$ is holomorphic in the
parameters of $W_1, W_2$.
\end{coro}
Summarizing,
\begin{prop}\label{specialization} Let the data be as
above. Let $S$ a finite collection of places containing the
archimedean place such that for $v \notin S$, the local data at
$v$ is unramified. Then we have
\begin{equation*}\begin{split}
\zZ(\varphi_1, \varphi_2, \mu \vertdot^s) = L(s, \pi_1, \mu)
& L(1-s, \pi_2, \mu^{-1}) \\
& \left\{\prod_v \Psi(W_1, W_2, f; \mu_v, \vertdot_v^{s},
\vertdot_v^{1-s})\right\}
\end{split}\end{equation*}
where by lemmas \ref{unramified-holo} and \ref{quotient-holo} the
expression in curly braces is a finite product and is entire.
\end{prop}
\section{The pull-back of the Whittaker function}
In this section, we aim to relate the local Euler factor of the
integral of Novodvorsky at the archimedean place to the
corresponding Euler factor of the integral considered in Section
\ref{Section:Integral}. For this purpose, we start by studying the
Whittaker function associated to $\theta (\varphi_1, \varphi_2;
f)$, and from that we derive formulae for the corresponding local
Whittaker functions.

\subsection{The Whittaker function} In this section we compute the Whittaker function of a
cuspidal function $\theta(\varphi_1, \varphi_2; f)$. Fix a
non-trivial character $\psi$ of $F \backslash \A$. Define a
character, again denoted by $\psi$, of the unipotent radical of
the Borel subgroup of $\GSp(4)$ by the following
\begin{equation*}
\psi( \begin{pmatrix} 1&v \\ &1 \\ && 1 \\ &&-v & 1
\end{pmatrix}
\begin{pmatrix} 1&& s & r \\ & 1 & r & t \\ &&1 \\ &&&1
\end{pmatrix}) = \psi(v + t).
\end{equation*}
Then we set
\begin{equation*}
W(g) = \int_{N(F) \backslash N(\A)} \theta(\varphi_1, \varphi_2;
f)(ng) \psi^{-1}(n) \, dn.
\end{equation*}
The $h^1$ and $h^2$ above can be taken to be $\begin{pmatrix} v(g)
\\ & 1 \end{pmatrix}$ and the identity matrix, respectively. Then
we have
\begin{thm}\label{Whittaker-theta} If $\tilde \pi_1 \ne \bar{\pi}_2$, we have
\begin{equation*}\begin{split}
W(g) = & \int_{\hat{N}(\A) \backslash H_1(\A)} W_1(\epsilon h_1
h^1) W_2(h_2 h^2) \\
& \omega(g, h_1 h^1, h_2 h^2) f (
\begin{pmatrix} 0 & -1 \\ 0& 0 \end{pmatrix}, I_{ 2 \times 2}) \,
dh_1 \, dh_2, \end{split}\end{equation*} where
\begin{equation*}
\hat{N} = \{ ( \begin{pmatrix} 1 & x \\ & 1 \end{pmatrix} ,
\begin{pmatrix} 1 & x \\ & 1 \end{pmatrix})\, \vert \, x \in \mathbb{G}_a
\}.
\end{equation*}
\end{thm}
\begin{proof} We start by
\begin{equation*}\begin{split}
W(g) = & \int_{H_1(F) \backslash H_1(\A)} \varphi_1(h_1
h^1)\varphi_2(h_2 h^2) \\
& \bigl( \sum_{M_1, M_2} \int_{N(F) \backslash N(\A)} \omega(ng;
h_1h^1, h_2h^2) f(M_1, M_2) \psi^{-1}(n) \, dn \bigr) \\
& d(h_1, h_2).
\end{split}\end{equation*}
Therefore, we have to study the expression
\begin{equation*}
I(M_1, M_2) =\int_{N(F) \backslash N(\A)} \omega(ng; h_1h^1,
h_2h^2) f(M_1, M_2) \psi^{-1}(n) \, dn.
\end{equation*}
For this we have
\begin{align*}
I(M_1, M_2)  =&  \int_{F \backslash \A} \bigl( \int_{(F \backslash
\A)^3} \omega (\begin{pmatrix} 1&& s & r \\ & 1 & r & t \\ &&1 \\
&&&1 \end{pmatrix}, I_2, I_2)  \\
& \omega(\begin{pmatrix} 1&v \\ &1 \\
&& 1 \\ &&-v & 1 \end{pmatrix} g, h_1 h^1 , h_2 h^2) f(M_1, M_2)
\psi^{-1}(t) \, dr \, ds \, dt \bigr) \\
&\psi^{-1}(v) \, dv \\
=& \int_{F \backslash \A}\omega(\begin{pmatrix} 1&v \\ &1 \\
&& 1 \\ &&-v & 1 \end{pmatrix} g, h_1 h^1 , h_2 h^2) f(M_1, M_2)
\\
& \bigl(\int_{(F \backslash \A)^3}\psi( tr \bigl( \begin{pmatrix}
s & r
\\ r & t \end{pmatrix}
\begin{pmatrix} \det M_1 & B(M_1, M_2) \\ B(M_2, M_1) & \det M_2 -1
\end{pmatrix} \bigr)
\, dr \, ds \, dt \bigr) \\
&\psi^{-1}(v) \, dv.
\end{align*}
But the inner most integral
\begin{equation*}
\int_{(F \backslash \A)^3}\psi( tr \bigl( \begin{pmatrix} s & r
\\ r & t \end{pmatrix}
\begin{pmatrix} \det M_1 & B(M_1, M_2) \\ B(M_2, M_1) & \det M_2 -1
\end{pmatrix} \bigr)
\, dr \, ds \, dt  = 0
\end{equation*}
unless $ \det M_1 =0$, $\det M_2 =1$, and $B(M_1, M_2) =0$, in
which case it is equal to $1$.
\begin{lem} Under the action of $H_1(F)$, the set $\mathcal{S}$
consisting of the pairs of matrices $(M_1, M_2)$ satisfying the
conditions just mentioned is the union of the following two
orbits:
\begin{enumerate}
\item The orbit of $(O, I)$. The stabilizer of this element is
the diagonal embedding of $\PGL(2)$ into $H_1$.
\item The orbit of $(\begin{pmatrix} 1 \\ & \end{pmatrix} ,
\begin{pmatrix} & 1 \\ -1 \end{pmatrix})$. The stabilizer of this
element is the subgroup $\tilde{N}$ of $H_1$ consisting of pairs
of matrices of the form
$$( \begin{pmatrix} 1 & x \\ & 1
\end{pmatrix}, w
\begin{pmatrix} 1 & x\\ 1 \end{pmatrix} w^{-1}).$$
\end{enumerate}
\end{lem}
\begin{proof} Since $\det M_1=0$, there are two cases to be
considered:
\begin{enumerate}
\item $M_1=0$,
\item $M_1 \ne 0$.
\end{enumerate}
It's obvious that the first case corresponds to the first orbit in
the statement of the lemma. Also the statement regarding the
stabilizer is immediate. Next we consider the case when $M_1 \ne
0$. It is clear that under the action of $H_1$, $M_1$ is
equivalent to the matrix $\begin{pmatrix} 1&0 \\ 0&0
\end{pmatrix}$. Next suppose $M_2 = \begin{pmatrix} a & b \\ c & d
\end{pmatrix}$. Since $B(M_1, M_2) = 0$ and $\det M_1 =0$, we
obtain that $\det (M_1 + M_2) =1$. This then implies that $d=0$.
But then since $\det M_2 = 1$, we obtain $c = - b^{-1}$. Hence
$M_2 = \begin{pmatrix} a & b \\ - b^{-1} \end{pmatrix}$. Next
consider the element
\begin{equation*}
h = ( \begin{pmatrix} 1 \\ & b^{-1} \end{pmatrix} \begin{pmatrix}
b^{-1} & a \\ & b \end{pmatrix}, \begin{pmatrix} b^{-1} \\ & 1
\end{pmatrix} ) \in H_1(F).
\end{equation*}
Then it is easy to check that
\begin{equation*}
h. ( \begin{pmatrix} 1 & 0 \\ 0 & 0 \end{pmatrix},\begin{pmatrix}
a & b \\ - b^{-1} \end{pmatrix}) = (\begin{pmatrix} 1 & 0 \\ 0 & 0
\end{pmatrix}, \begin{pmatrix}  & 1 \\ -1 \end{pmatrix}).
\end{equation*}
The statement regarding the stabilizer is straightforward.
\end{proof}
Next we study the contribution of each orbit to the Whittaker
integral. Corresponding to the two orbits obtained above, we have
the following two integrals:
\begin{equation*}\begin{split}
I_1 (g) =&  \int_{G(F) \backslash H_1(\A)} \int_{ F \backslash \A}
\omega( \begin{pmatrix} 1 & v \\ & 1 \\ && 1 \\ && -v & 1
\end{pmatrix}g, h_1 h^1 , h_2 h^2)\\
& f (\begin{pmatrix} 0 & 0 \\ 0 & 0 \end{pmatrix}, \begin{pmatrix}
1 \\ & 1 \end{pmatrix}) \varphi_1(h_1 h^1) \varphi_2(h_2 h^2)
\psi^{-1}(v) \, dv \, d(h_1, h_2),
\end{split}\end{equation*}
and
\begin{equation*}\begin{split}
I_2 (g) =&  \int_{\tilde{N}(F) \backslash H_1(\A)} \int_{ F
\backslash \A} \omega(
\begin{pmatrix} 1 & v \\ & 1 \\ && 1 \\ && -v & 1
\end{pmatrix}g, h_1 h^1 , h_2 h^2)\\
& f (\begin{pmatrix} 1 & 0 \\ 0 & 0 \end{pmatrix}, \begin{pmatrix}
&1 \\ -1 \end{pmatrix}) \varphi_1(h_1 h^1) \varphi_2(h_2 h^2)
\psi^{-1}(v) \, dv \, d(h_1, h_2).
\end{split}\end{equation*}
Then it is clear that
\begin{equation*}
W(g) = I_1(g) + I_2(g).
\end{equation*}
\begin{lem}
We have
\begin{equation*}
I_1(g ) =0,
\end{equation*}
except when $\tilde \pi_1 = \bar \pi_2$.
\end{lem}
\begin{proof}
By \cite{Harris-Kudla}, we have
\begin{align*}
I_1(g)=&  \int_{G(F) \backslash H_1(\A)} \int_{ F \backslash \A}
\omega( \begin{pmatrix} 1 & v \\ & 1 \\ && 1 \\ && -v & 1
\end{pmatrix}g \begin{pmatrix} I \\ & \nu(g)^{-1} I \end{pmatrix}) \\
& L( h_1 h^1 , h_2 h^2) f (\begin{pmatrix} 0 & 0 \\ 0 & 0
\end{pmatrix}, \begin{pmatrix} 1 \\ & 1 \end{pmatrix})
\varphi_1(h_1 h^1) \varphi_2(h_2 h^2) \\
& \psi^{-1}(v) \, dv \, d(h_1, h_2)
\end{align*}
\begin{align*}
=&  \int_{G(\A) \backslash H_1(\A)} \int_{\PGL_2(F) \backslash
\PGL_2(\A)}\int_{ F \backslash \A} \omega(
\begin{pmatrix} 1 & v \\ & 1 \\ && 1 \\ && -v & 1
\end{pmatrix}g \begin{pmatrix} I \\ & \nu(g)^{-1} I \end{pmatrix}) \\
& L( \gamma h_1 h^1 ,\gamma h_2 h^2) f (\begin{pmatrix} 0 & 0 \\ 0
& 0 \end{pmatrix}, \begin{pmatrix} 1 \\ & 1 \end{pmatrix})
\varphi_1(\gamma h_1 h^1) \varphi_2(\gamma h_2 h^2) \\
& \psi^{-1}(v) \, dv \, d\gamma \, d(h_1, h_2)
\end{align*}
\begin{align*}
=&\int_{G(\A) \backslash H_1(\A)} \int_{ F \backslash \A} \omega(
\begin{pmatrix} 1 & v \\ & 1 \\ && 1 \\ && -v & 1
\end{pmatrix}g \begin{pmatrix} I \\ & \nu(g)^{-1} I \end{pmatrix}) \\
& L(h_1 h^1 ,h_2 h^2) f (\begin{pmatrix} 0 & 0 \\ 0 & 0
\end{pmatrix}, \begin{pmatrix} 1 \\ & 1 \end{pmatrix}) \psi^{-1}(v)\\
&\left(\int_{\PGL_2(F) \backslash \PGL_2(\A)}\varphi_1(\gamma h_1
h^1) \varphi_2(\gamma h_2 h^2) \, d\gamma \right) \, dv  \, d(h_1,
h_2).
\end{align*}
The inner most integral
\begin{equation*}\begin{split}
\int_{\PGL_2(F) \backslash \PGL_2(\A)}& \varphi_1(\gamma h_1 h^1)
\varphi_2(\gamma h_2 h^2) \, d\gamma  = \\
& < \pi_1(h_1 h^1 ) \varphi_1 , \overline{\pi_2 (h_2 h^2)
\varphi_2}>_{L^2(\PGL_2(F) \backslash \PGL_2(\A))}.
\end{split} \end{equation*}
The statement of the lemma is now obvious.
\end{proof}

Next we study the contribution of the second orbit.
\begin{lem} We have
\begin{equation*}\begin{split}
I_2(g) = & \int_{\hat{N}(\A) \backslash H_1(\A)} W_{\varphi_1}(
\begin{pmatrix} 1 \\ & -1 \end{pmatrix} h_1
\begin{pmatrix} \nu(g) \\ & 1 \end{pmatrix}) W_{\varphi_2}(h_2) \\
& \omega(g, h_1 \begin{pmatrix} \nu(g) \\ & 1 \end{pmatrix}, h_2)
f (\begin{pmatrix} 0 & -1 \\ 0 & 0 \end{pmatrix}, I) \, d (h_1,
h_2).
\end{split}\end{equation*}
\end{lem}
In this lemma, $\hat{N}$ is the diagonal embedding of the
unipotent upper triangular matrices in $\GL(2)$ in $H_1$. Also if
$\varphi$ is a cuspidal automorphic function on $\GL_2(\A)$, we
have set
\begin{equation*}
W_\varphi(g) = \int_{ F \backslash \A} \varphi( \begin{pmatrix} 1 & x \\
& 1 \end{pmatrix} g ) \psi^{-1}(x) \, dx.
\end{equation*}
\begin{proof} The proof consists of simple
manipulations of the original expression for $I_2(g)$. We have
\begin{equation*}\begin{split}
I_2 (g) =&  \int_{\tilde{N}(F) \backslash H_1(\A)} \int_{ F
\backslash \A} \omega(
\begin{pmatrix} 1 & v \\ & 1 \\ && 1 \\ && -v & 1
\end{pmatrix}g, h_1 h^1 , h_2 h^2)\\
& f (\begin{pmatrix} 1 & 0 \\ 0 & 0 \end{pmatrix}, \begin{pmatrix}
&1 \\ -1 \end{pmatrix}) \varphi_1(h_1 h^1) \varphi_2(h_2 h^2)
\psi^{-1}(v) \, dv \, d(h_1, h_2).
\end{split}\end{equation*}
We recall that $\tilde{N}(F) = \{ ( \begin{pmatrix} 1 & x \\ & 1
\end{pmatrix}, w \begin{pmatrix} 1 & x \\ & 1 \end{pmatrix}
w^{-1}) \}$, and also that $h^1 = \begin{pmatrix} \nu(g) \\ & 1
\end{pmatrix}$ and $h^2 = I$. Using the formulae in
\cite{Harris-Kudla}, we have
\begin{equation*} \begin{split}
\omega(
\begin{pmatrix} 1 & v \\ & 1 \\ && 1 \\ && -v & 1
\end{pmatrix}g, & h_1 h^1 , h_2 h^2)
 f (\begin{pmatrix} 1 & 0 \\ 0 & 0 \end{pmatrix}, \begin{pmatrix}
&1 \\ -1 \end{pmatrix}) = \\
&\omega(g, h_1 h^1 , h_2 h^2)f (\begin{pmatrix} 1 & 0 \\ 0 & 0
\end{pmatrix}, \begin{pmatrix} 1 & -v \\ & 1 \end{pmatrix} \begin{pmatrix} & 1 \\ -1
\end{pmatrix}).
\end{split}\end{equation*}
Hence
\begin{align*}
I_2(g) = &\int_{\tilde{N}(F) \backslash H_1(\A)} \int_{ F
\backslash \A}\omega(g, h_1 h^1 , h_2 h^2)f (\begin{pmatrix} 1 & 0
\\ 0 & 0 \end{pmatrix}, \begin{pmatrix} 1 & -v \\ & 1 \end{pmatrix} \begin{pmatrix} & 1 \\ -1
\end{pmatrix}) \\
& \varphi_1(h_1 h^1)\varphi_2(h_2 h^2) \psi^{-1}(v) \, dv \,
d(h_1, h_2) \\
= &\int_{\tilde{N}(\A) \backslash H_1(\A)} \int_{ F \backslash \A}
\int_{ F \backslash \A}  \omega(g, \begin{pmatrix} 1 & u \\ & 1
\end{pmatrix} h_1 h^1 , w \begin{pmatrix} 1 & u \\ & 1 \end{pmatrix} w^{-1} h_2 h^2)\\
& f (\begin{pmatrix} 1 & 0 \\ 0 & 0 \end{pmatrix}, \begin{pmatrix}
1 & -v \\ & 1 \end{pmatrix} \begin{pmatrix} & 1 \\ -1
\end{pmatrix})
\varphi_1(\begin{pmatrix} 1 & u \\ & 1 \end{pmatrix} h_1
h^1)\varphi_2(w \begin{pmatrix} 1 & u \\ & 1 \end{pmatrix} w^{-1}
h_2 h^2)\\
&  \psi^{-1}(v) \, du \, dv \, d(h_1, h_2)
\end{align*}
Next by definition and Lemma 5.1.2 of \cite{Harris-Kudla}
\begin{equation*}
\omega(g,  \begin{pmatrix} 1 & u \\ & 1
\end{pmatrix} h_1 h^1 , w \begin{pmatrix} 1 & u \\ & 1 \end{pmatrix} w^{-1} h_2
h^2)
\end{equation*}
\begin{equation*}
= \omega ( g \begin{pmatrix} I \\ & \nu(g)^{-1} I
\end{pmatrix} )L (\begin{pmatrix} 1 & u \\ & 1
\end{pmatrix} h_1 h^1 , w \begin{pmatrix} 1 & u \\ & 1 \end{pmatrix} w^{-1} h_2
h^2)
\end{equation*}
\begin{equation*}
 =L (\begin{pmatrix} 1 & u \\ & 1
\end{pmatrix} h_1 h^1 , w \begin{pmatrix} 1 & u \\ & 1 \end{pmatrix} w^{-1} h_2
h^2)\omega ( \begin{pmatrix} I \\ & \nu(g)^{-1} I
\end{pmatrix} g ).
\end{equation*} This identity implies that
\begin{equation*}
\omega(g,  \begin{pmatrix} 1 & u \\ & 1
\end{pmatrix} h_1 h^1 , w \begin{pmatrix} 1 & u \\ & 1 \end{pmatrix} w^{-1} h_2 h^2)
f (\begin{pmatrix} 1 & 0 \\ 0 & 0 \end{pmatrix}, \begin{pmatrix} 1
& -v \\ & 1 \end{pmatrix} \begin{pmatrix} & 1 \\ -1 \end{pmatrix})
\end{equation*}
\begin{equation*}
= L(h_1 h^1 , h_2 h^2) \omega (\begin{pmatrix} I \\ & \nu(g)^{-1}
I \end{pmatrix} g ) f (\begin{pmatrix} 1 & 0 \\
0 & 0 \end{pmatrix}, \begin{pmatrix} 1 & -v \\ & 1 \end{pmatrix}
\begin{pmatrix} & 1 \\ -1 \end{pmatrix})
\end{equation*}
\begin{equation*}
=L(h_1 h^1 , \begin{pmatrix} 1 & -v \\ & 1 \end{pmatrix}
\begin{pmatrix} & -1 \\ 1 \end{pmatrix} h_2 h^2) \omega (\begin{pmatrix} I \\
& \nu(g)^{-1}
I \end{pmatrix} g ) f (\begin{pmatrix} 0 & -1 \\
0 & 0 \end{pmatrix}, I)
\end{equation*}
\begin{equation*}
 = \omega ( g,h_1 h^1 , \begin{pmatrix} 1
& -v \\ & 1
\end{pmatrix} \begin{pmatrix} & -1 \\ 1 \end{pmatrix} h_2 h^2)f (\begin{pmatrix} 0 & -1 \\
0 & 0 \end{pmatrix}, I).
\end{equation*}
Going back to $I_2(g)$, we obtain
\begin{equation*}\begin{split}
I_2(g) =& \int_{\tilde{N}(\A) \backslash H_1(\A)} \int_{ F
\backslash \A} \int_{ F \backslash \A} \omega ( g,h_1 h^1 ,
\begin{pmatrix} 1 & -v \\ & 1 \end{pmatrix}
w h_2 h^2)f (\begin{pmatrix} 0 & -1 \\
0 & 0 \end{pmatrix}, I) \\
& \varphi_1(\begin{pmatrix} 1 & u \\ & 1 \end{pmatrix} h_1
h^1)\varphi_2(w \begin{pmatrix} 1 & u
\\ & 1 \end{pmatrix} w^{-1} h_2 h^2) \psi^{-1}(v) \, du \, dv \, d(h_1, h_2)
\end{split}\end{equation*}
Next we make a change of variables $(h_1, h_2) \mapsto (h_1,
w^{-1} \begin{pmatrix} 1 & v \\ & 1 \end{pmatrix} h_2)$, to obtain
\begin{equation*}\begin{split}
I_2(g) =& \int_{\hat{N}(\A) \backslash H_1(\A)} \int_{ F
\backslash \A} \int_{ F \backslash \A} \omega ( g,h_1 h^1 , h_2 h^2)f (\begin{pmatrix} 0 & -1 \\
0 & 0 \end{pmatrix}, I) \\
& \varphi_1(\begin{pmatrix} 1 & u \\ & 1 \end{pmatrix} h_1
h^1)\varphi_2(\begin{pmatrix} 1 & u +v
\\ & 1 \end{pmatrix} h_2 h^2) \psi^{-1}(v) \, du \, dv \, d(h_1,
h_2).
\end{split}\end{equation*}
Now a change of variables $v \mapsto v-u$ and re-arranging the
order of integrals gives the result.
\end{proof}

Combining everything finishes the proof of the theorem.
\end{proof}
\subsection{Local Whittaker functions}\label{local-whittaker-functions}
In this paragraph, we study the integrals of the previous section
in some detail.

Suppose $\pi_1$ and $\pi_2$ are two irreducible admissible
representations of the group $\GL(2)$ over a local field, such
that $\tilde{\pi}_1 \ne \pi_2, \bar{\pi}_2$, and $\omega_{\pi_1}.
\omega_{\pi_2} =1$. For $W_i \in \mathcal{W}(\pi_i, \psi)$, for
$i=1, 2$, set
\begin{equation*}
\begin{split}
\mathbb{W}_v(W_1, W_2; f)(g) = & \int_{\hat{N}(F_v) \backslash
H_1(F_v)} W_1(\epsilon h_1
\begin{pmatrix} \nu(g) \\ & 1 \end{pmatrix}) W_2(h_2) \\
& \omega(g, h_1 \begin{pmatrix} \nu(g) \\ & 1 \end{pmatrix}, h_2)
f ( \begin{pmatrix} 0 & -1 \\ 0& 0 \end{pmatrix}, I_{ 2 \times 2})
\, dh_1 \, dh_2. \end{split}
\end{equation*}
\begin{prop} For all $W_i \in \mathcal{W}(\pi_i, \psi)$, $i=1,2$,
$K$-finite $f$ in the space of Schwartz-Bruhat functions, and $g
\in \GSp_4(F_v)$, the integral defining $\mathbb{W}(W_1, W_2;
f)(g)$ is absolutely convergent.
\end{prop}
\begin{proof} As usual we prove the proposition for the archimedean place.
It is clear that we only need to prove the absolute
convergence for $g = I_{4 \times 4}$. In order to do this, we
start by identifying a measurable set of representatives for
$\hat{N}(\R) \backslash H_1(\R)$, and identifying the
corresponding measure. On $H_1(\R)$, we have the following natural
set of representatives
\begin{equation*}
( \begin{pmatrix} 1 & x \\ & 1 \end{pmatrix} k_1,
\begin{pmatrix} 1 & y \\ & 1 \end{pmatrix}
\begin{pmatrix} \eta \\ & \eta^{-1} \end{pmatrix} k_2 ),
\end{equation*}
with $x,y \in \R$, $\epsilon \in \R^\times$, and $k_1, k_2 \in
\SO(2)$. Also the corresponding measure is
\begin{equation*}
\vert \eta \vert^{-2} \, dx \, dy \, d^\times \eta \, dk_1 \,
dk_2.
\end{equation*}
This statement implies that the set of elements of the form
\begin{equation*}
( \begin{pmatrix} 1 & x \\ & 1 \end{pmatrix} k_1,
\begin{pmatrix} \eta \\ & \eta^{-1} \end{pmatrix} k_2 ),
\end{equation*}
constitutes a measurable set of representatives for $\hat{N}(\R)
\backslash H_1(\R)$. Also with this normalization the measure is
\begin{equation*}
\vert \eta \vert^{-2} \, dx \, d^\times \eta \, dk_1 \, dk_2.
\end{equation*}
Hence we are reduced to proving the convergence of the following
integral:
\begin{equation*}\begin{split}
& \int_{K}\int_{K}\int_{\R}\int_{\R^\times} \left| W_1( \epsilon
\begin{pmatrix} 1 & x \\ & 1 \end{pmatrix} k_1)
W_2( \begin{pmatrix} \eta \\& \eta^{-1} \end{pmatrix}k_2) \right| . \\
& \left| f (k_1^{-1}\begin{pmatrix} 1 & - x \\ & 1 \end{pmatrix}
\begin{pmatrix} 0 & -1 \\ 0 & 0 \end{pmatrix}
\begin{pmatrix} \eta \\ & \eta^{-1} \end{pmatrix} k_2,
k_1^{-1}\begin{pmatrix} 1 & -x \\ & 1 \end{pmatrix}
\begin{pmatrix} \eta \\ & \eta^{-1} \end{pmatrix} k_2) \right| \\
& \, d^\times \eta \, dx \, dk_1 \, dk_2.
\end{split}\end{equation*}
Next we observe that in order to prove the absolute convergence of
this integral, we just need to prove the absolute convergence of
the integral over $\eta$ and $x$. Also since
\begin{equation*}
W_1( \epsilon
\begin{pmatrix} 1 & x \\ & 1 \end{pmatrix} k_1) = \psi(-x)
W_1(\epsilon k_1),
\end{equation*}
we obtain
\begin{equation*}
\left| W_1( \epsilon
\begin{pmatrix} 1 & x \\ & 1 \end{pmatrix} k_1) \right| = \left|
W_1(\epsilon k_1) \right|.
\end{equation*}
Hence we are reduced to proving the convergence of the following
integral:
\begin{equation*}\begin{split}
I= &\int_{\R}\int_{\R^\times} \left|
W_2( \begin{pmatrix} \eta \\& \eta^{-1} \end{pmatrix} )\right| . \\
& \left| f (\begin{pmatrix} 1 & - x \\ & 1 \end{pmatrix}
\begin{pmatrix} 0 & -1 \\ 0 & 0 \end{pmatrix}
\begin{pmatrix} \eta \\ & \eta^{-1} \end{pmatrix} ,
\begin{pmatrix} 1 & -x \\ & 1 \end{pmatrix}
\begin{pmatrix} \eta \\ & \eta^{-1} \end{pmatrix}) \right| \, d^\times \eta \, dx .
\end{split}\end{equation*}
But this integral is equal to
\begin{align*}
I = &\int_{\R}\int_{\R^\times} \left|
\omega_{\pi_2}(\eta^{-1})W_2(
\begin{pmatrix} \eta^2
\\& 1 \end{pmatrix} ) f (\begin{pmatrix} 0 & -
\eta^{-1} \\ 0 & 0 \end{pmatrix},
\begin{pmatrix} \eta & -x \eta^{-1} \\0 & \eta^{-1} \end{pmatrix}) \right|
\, d^\times \eta \, dx
\end{align*}
Now we write
\begin{equation*}
f (\begin{pmatrix} 0 & - \eta^{-1} \\ 0 & 0 \end{pmatrix},
\begin{pmatrix} \eta & -x \eta^{-1} \\0 & \eta^{-1} \end{pmatrix}) =
q(\eta, \eta^{-1}, x \eta^{-1}),
\end{equation*}
where $q$ is some Schwartz-Bruhat function in three variables. We
then need to prove the convergence of the integral
\begin{align*}
I = &\int_{\R}\int_{\R^\times} \left|
\omega_{\pi_2}(\eta^{-1})W_2(
\begin{pmatrix} \eta^2
\\& 1 \end{pmatrix} ) q(\eta, \eta^{-1}, x \eta^{-1}) \right|
\, d^\times \eta \, dx,
\end{align*}
which after a change of variables $x \mapsto x \eta$ and
integration over $x$ is equivalent to the convergence of an
integral of the form
\begin{equation*}
\int_{\R_+^\times} \left| W( \begin{pmatrix} \eta \\ & 1
\end{pmatrix} \xi(\eta) \right| \eta^\sigma \, d ^\times \eta
\end{equation*}
for $\xi \in \mathcal{S}(\R^\times)$. Such an integral always
converges by the moderate growth of the Whittaker function.
\end{proof}

Going back to the global situation, we choose $\varphi_i$, for
$i=1,2$, so that
\begin{equation*}
W_{\varphi_i} = \otimes_{v \in \mathcal{M}_F} W_v^i.
\end{equation*}
We also choose $f$ to be a pure tensor of the form $\otimes_v
f_v$. Then theorem \ref{Whittaker-theta} can be written in the
form
\begin{equation*}
W(g) = \prod_v \mathbb{W}_v(W_v^1, W_v^2; f_v)(g_v).
\end{equation*}
under appropriate conditions. This implies that for each local
place $v$, if $W_v$ is a $K_v$-finite vector in the local
Whittaker model, there is a choice of the data such that $W_v
=\mathbb{W}_v(W_v^1, W_v^2; f_v)$. It is clear from the
construction that, in the archimedean situation, the space of all
such ${\mathbb W}$'s forms a $(\mathfrak{g}, K)$-module.

\section{Archimedean Zeta function} In this section, we use the results of the previous
paragraphs to obtain information about the archimedean zeta
function. We have by lemma \ref{Bessel-Zeta}
\begin{equation}\begin{split}\label{Bessel-Zeta-repeat}
B(\phi, \chi_s)=\int_{{\mathbb A}^\times} \int_{\mathbb A}
W_\phi\left(
\begin{pmatrix} y \\ &y \\ &&1 \\ &x&&1 \end{pmatrix}w^{-1}\right) \mu(y)
\vert y \vert ^{s - {3 \over 2}} \, dx \, d^\times y,
\end{split}\end{equation}
with
\begin{equation*}
w = \begin{pmatrix} 1 \\ &&&1 \\ &&1 \\ &-1 \end{pmatrix}
\end{equation*}
and
\begin{equation*}
\chi_s(y) = \mu(y)|y|^{s - \frac{1}{2}}.
\end{equation*}
If we set $\phi = \theta(\varphi_1, \varphi_2; f)$, the left hand
side of the above identity will be equal to what we have called
$\zZ(\varphi_1, \varphi_2, f; \mu \vertdot^s)$. We saw in
\ref{specialization} that
\begin{equation*}\begin{split}
\zZ(\varphi_1, \varphi_2, \mu \vertdot^s) = L(s, \pi_1, \mu)
& L(1-s, \pi_2, \mu^{-1}) \\
& \left\{\prod_v \Psi(W^v_1, W^v_2, f; \mu_v \vertdot_v^{s},
\mu_v^{-1} \vertdot_v^{1-s})\right\}.
\end{split}\end{equation*}
If we choose our vectors appropriately, that is factorizable, the
right hand side of (\ref{Bessel-Zeta-repeat}) is now equal to
\begin{align*}
\prod_v & Z_{v, N}(s, \pi_v(w^{-1})\mathbb{W}_v(W_v^1, W_v^2;
f_v),\mu_\infty) \\
& = Z_{\infty, N}(s,
\pi_\infty(w^{-1})\mathbb{W}_\infty(W_\infty^1, W_\infty^2;
f_\infty), \mu_\infty)\\
& \,\,\,\,\,\,\,\,\times L_S(s, \Pi, \mu) \prod_{v \in S
\backslash \{ \infty \}} Z_{v, N}(s,
\pi_v(w^{-1})\mathbb{W}_v(W_v^1, W_v^2; f_v), \mu_v)
\end{align*}
By the main result of \cite{Takloo-Bighash}, for each local place
$v \in S \backslash \{\infty\}$, we can choose $W^{\tt sp}_v \in
\mathcal{W}(\Pi_v)$ in such a way that
\begin{equation*}
Z_{v, N}(s, \Pi_v^{-1}(w^{-1}) W^{\tt sp}_v , \mu_v) = L_v (s,
\Pi_v, \mu_v).
\end{equation*}
By the remark at the end of \ref{local-whittaker-functions}, we
can choose the local data such that
\begin{equation*}
\mathbb{W}_v(W_v^1, W_v^2; f_v) = W^{\tt sp}_v.
\end{equation*}
With this choice of the local data, we have
\begin{equation}\begin{split}\label{local-global}
& Z_{\infty, N}(s, \pi_\infty(w^{-1})\mathbb{W}_\infty(W_\infty^1,
W_\infty^2; f_\infty), \mu_\infty)
 \\
& = \Phi^{\tt finite}_S(\pi_1, \pi_2, \mu, s; W_1, W_2, f)
L_\infty(s, \pi_1,
\mu) L_\infty(1-s, \pi_2, \mu^{-1})\\
& \times \Psi(W^\infty_1, W^\infty_2, f_\infty; \mu_\infty
\vertdot_\infty^{s}, \mu_\infty^{-1} \vertdot_\infty^{1-s}),
\end{split}\end{equation}
with
\begin{equation*}\begin{split}
\Phi_S^{\tt finite}& (\pi_1, \pi_2, \mu, s; W_1, W_2, f) \\
& =\frac{L^\infty(s, \pi_1, \mu) L^\infty(1-s, \pi_2,
\mu^{-1})}{L^\infty(s, \Pi, \mu)} \prod_{v \in S \backslash \{
\infty \}} \Psi(W^v_1,
W^v_2, f; \mu_v \vertdot_v^{s}, \mu_v^{-1} \vertdot_v^{1-s}) \\
& =\prod_{v \in S \backslash \{ \infty \}} \Psi(W^v_1, W^v_2, f;
\mu_v \vertdot_v^{s}, \mu_v^{-1} \vertdot_v^{1-s}),
\end{split}\end{equation*}
if $\mu$ is chosen in such a way that for $v \in S \backslash \{
\infty \}$, the local quasi-character $\mu_v$ is highly ramified.
Combining everything proves the first statement of the following
theorem:
\begin{thm}\label{ratio-arch}
In the above situation, for each $K$-finite $W \in
\mathcal{W}(\Pi_\infty)$, the ratio
\begin{equation*}
\frac{Z(s, W, \mu_\infty)}{L_\infty(s, \pi_1^\infty, \mu_\infty)
L_\infty(s, \tilde{\pi}_2^\infty, \mu_\infty)}
\end{equation*}
extends to an entire function of $s$. Furthermore, for each $s$,
there is a choice of $W$ such that the above expression does not
vanish at $s$.
\end{thm}
\begin{proof}
We only need to prove the second statement. In order to do this,
we prove the existence of an entire function $\Phi(s)$ such that
\begin{equation}\begin{split}
Z_{\infty, N}&(s, \pi_\infty(w^{-1})\mathbb{W}_\infty(W_\infty^1,
W_\infty^2; f_\infty), \mu_\infty)
 \\
=& \frac{1}{\Phi(s)} L_\infty(s, \pi_1,
\mu) L_\infty(1-s, \pi_2, \mu^{-1})\\
& \times \Psi(W^\infty_1, W^\infty_2, f_\infty; \mu_\infty
\vertdot_\infty^{s}, \mu_\infty^{-1} \vertdot_\infty^{1-s}).
\end{split}\end{equation}
By proposition \ref{quotient-holo} there is a choice of the data
with the required property. Again we assume that $\mu$ is highly
ramified for non-archimedean $v \in S$, and unramified outside
$S$. In order to show the existence of $\Phi(s)$ it is not hard to
see that if we can show the existence of local non-archimedean
data with the property that
\begin{equation*}
L_v(s, \pi_1, \mu)  L_v(1-s, \pi_2, \mu^{-1}) \Psi(W^v_1, W^v_2,
f; \mu_v \vertdot_v^{s}, \mu_v^{-1} \vertdot_v^{1-s})
\end{equation*}
is a constant, then we can take
\begin{equation*}
\Phi(s) = C \prod_{v \in S \backslash \{\infty \}}Z_{v, N}(s,
\pi_v(w^{-1})\mathbb{W}_v(W_v^1, W_v^2; f_v), \mu_v),
\end{equation*}
with $C$ the obvious non-zero constant. The existence of such data
is the statement of Corollary \ref{coro:constant}.

 We claim that the function $\Phi(s)$ is nowhere
vanishing. To see this, we set
\begin{equation*}
F_1(W_\infty^1, W_\infty^2, s) = \frac{Z_{\infty, N}(s,
\pi_\infty(w^{-1})\mathbb{W}_\infty(W_\infty^1, W_\infty^2;
f_\infty), \mu_\infty)}{L_\infty(s, \pi_1, \mu) L_\infty(1-s,
\pi_2, \mu^{-1})}
\end{equation*}
\begin{equation*}
F_2(W_\infty^1, W_\infty^2, s)= \Psi(W^\infty_1, W^\infty_2,
f_\infty; \mu_\infty \vertdot_\infty^{s}, \mu_\infty^{-1}
\vertdot_\infty^{1-s}).
\end{equation*}
So far we know that given any $W_\infty^1, W_\infty^2$, the
complex functions $F_1(s), F_2(s)$ are both entire. Next, let
$s_0$ be given and suppose $\Phi(s_0) = 0$; but,
\begin{equation}\label{phi-s}
F_2(s) = \Phi(s) F_1(s),
\end{equation}
which would then imply that for all choices of data we must have
$F_2(s_0) = 0$ which, by proposition \ref{quotient-holo}, is not
true. This finishes the proof of the theorem. \end{proof}
\begin{rem}\label{phi-ramified-S}
We observe that the function $\Phi(s)$ defined in the proof of the
theorem does not depend on $W_\infty^1, W_\infty^2$, and its
dependence on $\pi_1^\infty, \pi_2^\infty$ is merely through the
non-archimedean components of the automorphic representations
$\pi_1, \pi_2$. As $\Phi(s)$ is the product of polynomials of
$q_v^{-s}$, for $v \in S$, and as it nowhere vanishing, it is a
function of the form
$$
A B ^ {-s}
$$
with $B$ rational. Also prime numbers appearing in the
decomposition of $B$ are all from the set $S$. We will see later
that $\Phi(s)$ is in fact a constant.
\end{rem}
\subsection{Analytic continuation} Let $\tau$ be a complex number with $\Re \tau >0$.
Then one can consider the archimedean principal series
representation $\pi(\tau)= \Ind (\vertdot^{\tau} \otimes
\vertdot^{-\tau})$. Let $\rho_\tau: W_\R \to \GL_2(\C)$ be the $L$
parameter associated with the representation $\pi(\tau)$. We
observe that if $\pi(\tau)$ is irreducible, the corresponding $L$
packet has a single element. Then as described in
\cite{Casselman-Shalika} one can consider a continuous map
$$
P(\tau) : \mathcal{S}(\GL_2(\R)) \longrightarrow \pi(\tau).
$$
Also for $v \in \pi(\tau)$, we set
$$
W(v, g) = \int_{N(\R)} v(ng) \psi^{-1}(n) \, dn
$$
when the integral converges. Fix a Schwartz function $f$, and set
$$
W_\tau(f; g) := W( P(\tau)(f), g).
$$
A theorem of Shahidi asserts that $W_\tau$ extends to an entire
function of $\tau$. Usually, suppressing $f$, we simply write
$W_\tau$. Fix two sections of $W_\tau$, say $W_{\tau_1}$ and
$W_{\tau_2}$. Next, consider the function
$$
\mathbb{W}_f(\tau_1, \tau_2) := \mathbb{W}(W_{\tau_1}, W_{\tau_2};
f)
$$
as before. We write $F_i(\tau_1, \tau_2, s)$, $i=1, 2$, instead of
the functions of the previous paragraph.

Let $\C_{\rm aut}$ be the collection of those complex numbers
$\tau$ with the property that $\pi(\tau)$ occurs as the
archimedean component of some automorphic cuspidal representation
of the group $\GL(2)$. It is well-known that $\C_{\rm
temp}:=\C_{\rm aut} \cap i \R$ is dense in $i \R$.

The function $\mathbb{W}_f(\tau_1, \tau_2)$ is entire on $\C^2$,
and for fixed $(\tau_1, \tau_2) \in \C^2$ defines a Whittaker
function on $\GSp(4, \R)$. Also by construction if $\tau_1, \tau_2
\in \C_{\rm temp}$, the function $\mathbb{W}_f(\tau_1, \tau_2)$
will make up the $K$-finite Whittaker model of the unique element
of the local $L$ packet $\varphi(\rho_{\tau_1}, \rho_{\tau_2})$.
In fact, if we stay away from the points of reducibility, the
unique element of the $L$ packet given by $\varphi(\rho_{\tau_1},
\rho_{\tau_2})$ is generic.

We have established the identity
$$
F_1(\tau_1, \tau_2, s) = \Phi(s) F_2(\tau_1, \tau_2, s)
$$
whenever $(\tau_1, \tau_2) \in \C_{\rm aut} \times \C_{\rm aut}$,
and $\Re s > b(\tau_1, \tau_2)$. Presumably, the function
$\Phi(s)$ depends on $s$, and, though we have suppressed the
dependence, on $\tau_1, \tau_2$. We now show that for $\tau_1,
\tau_2 \in \C_{\rm temp}$, $\Phi(s)$ is an absolute constant
independent of all variables. For this we follow the argument of
lemma 5 of \cite{Venkatesh}, which is in the spirit of
Burger-Li-Sarnak. The proof of Lemma 5 of \cite{Venkatesh} implies
that given $\tau \in i\R$ one can find an automorphic cuspidal
representation of $\GL(2)$ with archimedean component arbitrarily
close to $\pi(\tau)$ and ramified only at one prescribed place.
This, applied to a pair of tempered representations of $\GL(2)$
considered as a representation of $\GO(2,2)$, implies that given a
tempered representation of $\GO(2,2)(\R)$ one can find two
automorphic cuspidal representations with disjoint sets $S$. This
observation combined with remark \ref{phi-ramified-S} proves that
$\Phi(s)$ must be a constant. Next, we have
$$
F_1(\tau_1, \tau_2, s) = \Phi F_2(\tau_1, \tau_2, s)
$$
whenever $\tau_1, \tau_2 \in \C_{\rm temp}$ and $\Re s > b(s_1,
s_2)$. The density of $\C_{\rm temp}$ in $i \R$ then implies that
the identity must hold for all $\tau_1, \tau_2$, whenever $\Re s >
b(\tau_1, \tau_2)$. But we have seen that $F_2$ is entire as a
function of three complex variables; consequently, as $F_1$ and
$F_2$ agree on an open set, $F_2$ is the analytic continuation of
$F_1$. Consequently, whatever we proved about $F_2$ carries over
to $F_1$.

{\em Address: Department of Mathematics, Princeton University,
Princeton, NJ 08544. \\ Email:} {\tt rtakloo@math.princeton.edu}
\end{document}